\documentclass[11pt]{article}

\usepackage[margin=1in]{geometry}
\usepackage{amsmath, amssymb, amsthm}
\usepackage{adjustbox}
\usepackage{mathrsfs}
\usepackage{graphicx}
\usepackage{booktabs}
\usepackage{enumitem}
\usepackage{bbm}
\usepackage{natbib}
\usepackage[colorlinks=true,citecolor=blue,urlcolor=blue,linkcolor=blue]{hyperref}
\def\var{{\mathrm{var}}}

\newcommand\independent{\protect\mathpalette{\protect\independenT}{\perp}}
\def\independenT#1#2{\mathrel{\rlap{$#1#2$}\mkern2mu{#1#2}}}

\theoremstyle{definition}
\newtheorem{definition}{Definition}

\theoremstyle{remark}

\newcommand{\ran}{\mathrm{ran}}

\theoremstyle{plain}

\newtheorem{theorem}{{\bf Theorem}}
\newtheorem{lemma}{{\bf Lemma}}
\newtheorem{assumption}{{\bf Assumption}}
\newtheorem{corollary}{{\bf Corollary}}
\newtheorem{theoremC}{{\bf Theorem}}
\newtheorem{proposition}[theoremC]{{\bf Proposition}}
\def\hii#1{\hi{(#1)}}

\def\one{\mathbbm{1}}

\def\boe{\begin{enumerate}}
\def\eoe{\end{enumerate}}
\def\dom{\mathrm{dom}}
\newcommand\ca[1]{{\cal{#1}}}
\newcommand\lo[1]{_{\nano{#1}}}
\newcommand\hi[1]{^{\nano{#1}}}

\def\L{{\cal L}}

\def\tsum{\textstyle{\sum}}

\def\trans{^{\mbox{\tiny{\sf  T}}}}

\def\inv{^{\mbox{\tiny $-1$}}}

\def\var{\mathrm{var}}

\newcommand{\indep}{\;\, \rule[0em]{.03em}{.65em} \hspace{-.41em}
\rule[-.02em]{.65em}{.03em} \hspace{-.41em}
\rule[0em]{.03em}{.65em}\;\,}

\def\tr{\mathrm{tr}}
\def\vec{\mathrm{vec}}
\def\mat{\mathrm{mat}}

\def\trans{^{\mbox{\tiny{\sf T}}}}

\def\ali{&\,}

\def\ran{\mathrm{ran}}

\def\of{{\nano {\circ}}}
\def\nano{\scriptscriptstyle}

\def\real{{\mathbb R}}

\def\ka{\kappa}

\def\oc{^{\nano \perp}}

\def\L2T{L \lo 2 (T)}
\def\L2TX{L \lo 2 (T\lo X)}
\def\L2TX{L \lo 2 (T\lo Y)}
\def\tsum{\textstyle{\sum}}
\def\ali{&\,}
\def\spn{{\rm{span}}}
\def\ali{&\,}

\def\eod{\end{document}}

\def\cran{{\overline{{\rm{ran}}}}}
\def\dom{\mathrm{dom}}

\newcommand{\red}[1]{{\leavevmode\color{red}{#1}}}

\newfont{\rsfsten}{rsfs10 scaled 1100}
\newfont{\rsfstena}{rsfs10 scaled 800}
\newfont{\rsfstenb}{rsfs10 scaled 500}

\def\dom{{\mathrm{dom}}}

\title{Structure-Preserving Nonlinear Sufficient Dimension Reduction for Tensors}
\author{
Dianjun Lin, 
Bing Li, and
Lingzhou Xue \\ Department of Statistics, The Pennsylvania State University
}
\date{} 

\begin{document}
\maketitle

\begin{abstract}
We introduce two nonlinear sufficient dimension reduction methods for regressions with tensor-valued predictors. Our goal is two-fold: the first is to preserve the tensor structure when performing dimension reduction, particularly the meaning of the tensor modes, for improved interpretation; the second is to substantially reduce the number of parameters in dimension reduction, and thereby achieving model parsimony and enhancing estimation accuracy. Our two tensor dimension reduction methods echo the two commonly used tensor decomposition mechanisms: one is the Tucker decomposition, which reduces a larger tensor to a smaller one; the other is the CP-decomposition, which represents an arbitrary tensor as a sequence of rank-one tensors. We developed the Fisher consistency of our methods at the population level and established their consistency and convergence rates. Both methods are easy to implement numerically: the Tucker-form can be implemented through a sequence of least-squares steps, and the CP-form can be implemented through a sequence of singular value decompositions. We investigated the finite-sample performance of our methods and showed substantial improvement in accuracy over existing methods in simulations and two data applications. 
\end{abstract}


\section{Introduction}

Sufficient dimension reduction (SDR) provides a statistical framework for high-dimensional regression by reducing the predictor dimension while retaining all information about the response; see, for example, \cite{Li1991, cook1991, li2018sufficient}. 
Classical linear SDR methods are primarily developed for vector-valued predictors; they were then extended to accommodate tensor-valued predictors. In this paper, we develop the nonlinear SDR for tensors to enhance estimation accuracy and achieve further reduction of dimensionality. 

Tensors, or multidimensional arrays \citep{kolda2009tensor}, are prevalent in modern real-world applications such as neuroimaging, computer vision, and natural language processing. 
A straightforward approach for applying SDR to tensor predictors is to vectorize them and then apply the conventional SDR techniques \citep{li2018sufficient}; however, this destroys the intrinsic tensor structure and leads to unnecessarily high dimensionality. 
To address this limitation, efforts have been devoted to developing SDR methods specifically tailored for matrix- and tensor-valued data; see \cite{Li2010, zhong2015, fan2017sufficient, liu2018, zeng2013, luo2022inverse, yu2022nonparametric, xing2024, zhang2024dimension}, and others. 

\cite{Li2010} was the first in the literature to incorporate tensor predictors into SDR by introducing the \emph{dimension folding} framework, which preserves the tensor structure and drastically reduces the number of parameters. 
Subsequent works extended this idea along various directions. 
\cite{pfeiffer2012} proposed an inverse dimension folding method for longitudinal data. 
\cite{ding2014} developed two inverse dimension folding approaches based on principal components and principal fitted components (DF-PCA and DF-PFC), followed by a tensor-sliced inverse regression method \citep{ding2015}. 
\cite{xue2015} proposed forward regression–based dimension folding techniques, and \cite{wang2022likelihood} introduced a likelihood-based estimator for tensor dimension folding. 
Other related approaches for tensor-valued data include matrix logistic regression \citep{hung2013} and tensor regression \citep{zhou2013}. 
More recently, \cite{2025_kapla_tensor_arxiv} proposed a generalized multilinear model framework for SDR on tensor-valued predictors, unifying several existing linear tensor SDR methods within a broader multilinear structure.
All these methods are \emph{linear} tensor SDR approaches, in which the sufficient predictors are linear functions of the tensor-valued original predictors.  

Another important line of research underlying the present development is \emph{nonlinear sufficient dimension reduction}, which addresses the conditional independence problem
\[
Y \independent X \mid \{ h_1(X), h_2(X), \ldots, h_d(X) \},
\]
where $X \in \mathbb{R}^p$ is the predictor, $Y$ is the response, and $h_1, \ldots, h_d$ are nonlinear functions belonging to an appropriate functional space. 
\cite{wu2008} and \cite{yeh2008} extended sliced inverse regression (SIR) to the nonlinear setting by introducing the kernel SIR (KSIR), which represents the sufficient predictors $h_1, \ldots, h_d$ through kernel functions. 
\cite{li2011psvm} proposed the principal support vector machine as another nonlinear SDR framework. 
\cite{Lee2013} established a unified theoretical foundation for nonlinear SDR based on the central $\sigma$-field, and formulated the generalized SIR (GSIR) and generalized sliced average variance estimator (GSAVE). 
\cite{li2017} extended GSIR and GSAVE to functional data, leading to f-GSIR and f-GSAVE. 
Recently, \cite{song2023} introduced a nonlinear extension of the principal fitted component method based on the likelihood of the predictors. For a comprehensive treatment of nonlinear SDR, see Chapters 12–14 of \cite{li2018sufficient}.

To enhance the flexibility and accuracy of dimension reduction for tensor-valued data, we develop two \emph{nonlinear tensor sufficient dimension reduction} (NTSDR) methods that combine the strengths of tensor structures and nonlinear representations. 
Our framework is constructed in the tensor product space of reproducing kernel Hilbert spaces (RKHSs) defined on the left and right singular vectors of matrices, or more generally, on the mode vectors of higher-order tensors. 
Following \cite{Li2010}, we extend the tensor envelope concept to the nonlinear setting, thereby linking nonlinear SDR for vector-valued and tensor-valued predictors.

The proposed methods correspond to the two classical tensor decompositions. 
The \emph{Tucker form}, inspired by the Tucker decomposition \citep{tucker1966three, hitchcock1927tensor}, generalizes the dimension folding method of \cite{Li2010} to the nonlinear case. 
The \emph{CP form} (for CANDECOMP/PARAFAC; see, e.g., \cite{hitchcock1927tensor, zhang2001rankone}) extends the linear tensor SDR framework of \cite{zhong2015}. 
For both methods, we establish the Fisher consistency at the population level and, using the CP-version as a prototype, we prove the consistency and convergence rates of the estimated sufficient predictors. 
At the sample level, the Tucker form is implemented by iterative least-squares updates with closed-form solutions, whereas the CP form relies on iterative singular value decompositions. 
Although our exposition focuses on matrix-valued predictors, the CP form naturally extends to higher-order tensors. 
Simulation studies and real applications show that the proposed methods substantially outperform existing SDR approaches.

The rest of the paper is organized as follows.
Section 2 introduces the population-level structures of the proposed methods.
Section 3 defines the tensor regression operator and the stretched-out GSIR. 
Sections 4 and 5 present the Tucker- and CP-form nonlinear tensor SDR methods. 
Section 6 develops the asymptotic theory for our estimators.
Section 7 extends the CP-form method from matrices to general tensors, and Section 8 details the sample-level implementation. 
Section 9 further generalizes the framework to high-order tensor structures and discusses computational scalability.
Finally, Sections 10 and 11 present simulation studies and a real data application, respectively. The detailed proofs are presented in the supplement.  

\def\hi{^}
\def\trans{\hi {\mathsf T}}
\def\lo{_}
\def\real{\mathbb{R}}
\def\indep{\independent}
\def\realLR{\real \hi {p  \times q}}

\section{Outline of our approaches}

This section proposes two approaches to nonlinear SDR for tensors: the Tucker approach and the CP approach.  
To motivate our nonlinear constructions, we begin with a brief outline of the linear SDR approaches.

\subsection{Linear SDR for tensors}
Let $X \in \mathbb{R}^{p \times q}$ be a matrix-valued predictor and $Y \in \mathbb{R}$ a scalar response. 
A straightforward approach to dimension reduction treats $X$ as the vector $\mathrm{vec}(X)$ and applies the standard linear SDR method \citep{li2018sufficient}, leading to
\begin{align}\label{eq:linear stretch}
    Y \,\independent\, \mathrm{vec}(X) \mid \eta^\top \mathrm{vec}(X),
\end{align}
where $\eta$ is a matrix in $\real \hi {pq\times d}$ with the target dimension \( d \) after reduction. As discussed in \cite{Li2010}, this vectorization approach involves an excessive number of parameters ($p \times q  \times  d $ in total) and the sufficient predictors do not have a good tensor interpretation, because the row and column effects are totally ignored. 
To address this issue, \cite{Li2010} proposed to replace (\ref{eq:linear stretch}) by 
\begin{align}\label{eq:tucker linear}
 Y \independent X \mid \alpha\trans  X \beta,   
\end{align}
where  \( \alpha \in \mathbb{R}^{p \times s } \), \( \beta \in \mathbb{R}^{q \times t} \) and $s \le p$, $t \le q$. In this formulation, the number of parameters is $p s + q t$, which is much smaller than $p q d$. Also, $\alpha$ corresponds to the row mode and $\beta$ to the column mode of the tensor. Hence, dimension reduction is achieved for each mode of the tensor. Note that (\ref{eq:tucker linear}) can be regarded as a structured version of (\ref{eq:linear stretch}), and it is equivalent to 
\[
Y \independent \vec( X) \mid (\beta \otimes \alpha)\trans  \vec(X).
\]  
\cite{zhong2015} proposed an alternative approach  to find  $\tau \lo 1 , \ldots, \tau \lo d  \in \real \hi {p}$ and $\xi \lo 1  , \ldots, \xi \lo d  \in \real \hi {q}$ satisfying   
\begin{align}\label{eq:linear cp}
    Y \independent \vec(X) \mid  \, \tau \lo 1  \trans X \xi \lo 1 , \ldots, \tau \lo d  \trans X \xi \lo d.
\end{align} 
Similar to the formulation of \cite{Li2010}, SDR in this setting is also carried out mode by mode,  which results in a substantial reduction in the number of parameters compared with \eqref{eq:linear stretch}.

In the context of tensor decomposition (see, e.g., \cite{hitchcock1927tensor} and \cite{kolda2009tensor}), problem \eqref{eq:tucker linear} resembles the Tucker decomposition, where a high-dimensional tensor is represented by a low-dimensional core tensor. In contrast, problem \eqref{eq:linear cp} corresponds to the rank-one or CP (canonical polyadic) decomposition. Thus, we refer to \eqref{eq:tucker linear} as the Tucker-form linear SDR for tensors and \eqref{eq:linear cp} as the CP-form linear SDR for tensors. Our nonlinear extensions will follow the same formulations.

\subsection{Nonlinear feature extraction} 
A critical step in many nonlinear dimension reduction methods is feature mapping, in which the original data are mapped into a higher-dimensional space. 
Linear dimension reduction is then performed in this space, and the result is subsequently projected back to the original space, yielding nonlinear dimension reduction. 
See, for example, \cite{Scholkopf1997,fukumizu2007statistical,li2018sufficient,zhang2024nonlinear}. 
We adopt this general strategy to extend linear SDR for tensors to nonlinear SDR for tensors, while preserving the intrinsic tensor structure during the nonlinear mapping. To this end, we first perform singular value decomposition on  $x \in \real\hi {p \times q} $: 
\begin{align*}
  x= \sum_{i=1}^r \lambda_i (x)  U_i \hi 0 (x)    V_i \hi 0  (x)  \trans ,  
\end{align*}
where $U \lo 1 \hi 0  (x) , \ldots, U \lo r \hi 0  (x)  \in \real \hi {p }$ are the left singular vectors  of $x$, $V \lo 1 \hi 0  (x) , \ldots, V \lo r \hi 0  (x) \in \real \hi {q}$ are right singular vectors of $x$, $\lambda \lo 1 (x), \ldots, \lambda \lo r(x) $ are the singular values of $x$, and   \( r = \min\{p,q\} \) is the rank of \( x\). Note that the above singular value decomposition is not unique: if we replace $U_i^0(x)$ with $-U_i^0(x)$ and $V_i^0(x)$ with $-V_i^0(x)$, the decomposition remains valid. 
To ensure the uniqueness of our construction and to avoid ambiguity in subsequent developments, we assume that the first component  of $U_i^0(x)$ is positive and the sign of the first component of $V_i^0(x)$ is such that $\lambda \lo i (x)$ is positive. 
This assumption may fail when the first components of $U_i^0(x)$ and/or $V_i^0(x)$ are zero, but this issue is negligible, since for any $X$ whose distribution dominates the Lebesgue measure, the probability that $U_i^0(X)$ or $V_i^0(X)$ has a zero first component is zero. 
Another condition that guarantees uniqueness is that all nonzero singular values of $X$ are distinct, which is also reasonable because, if $X$ is a continuous random matrix whose distribution dominates the Lebesgue measure, then its nonzero singular values are distinct with probability one. 
An alternative approach to ensuring uniqueness involves the concept of odd and even kernels; see \cite{virta2024}. 

The assumptions on the random matrix $X$ are summarized below; these are automatically satisfied if the distribution of $X$ dominates the Lebesgue measure. Under these conditions, the quantities $\lambda \lo i(X)$, $U \lo i \hi 0 (X)$, $V \lo i \hi 0 (X)$ are uniquely determined by $X$ with probability one.
\vspace{-.1in}

\begin{assumption}\label{ass:assumptions2}
The random matrix $X$ satisfies the following conditions: 
\begin{enumerate}
    \item the rank of $X$ is $r = \min (p,q)$ with probability 1; 
    \item the first components of $U \lo i(X)$ and $V \lo i(X)$, $i = 1, \ldots, r$,  are nonzero with probability 1; 
    \item the nonzero singular values $\lambda \lo 1 (X), \ldots, \lambda \lo r (X)$ are simple with probability 1; that is, the left- or right-singular vector spaces corresponding to $\lambda \lo i (X)$, $i = 1, \ldots, r$,  have dimension 1 with probability 1. 
\end{enumerate}
\end{assumption} \vspace{-.1in}

Next, we define the reproducing kernel Hilbert spaces (RKHSs) for the left and right singular vectors, as well as their tensor product space.  Let \(    \ka \lo U : \real \hi {p  } \times 
    \real \hi {p  } \to \real,
    \ka \lo V : \real \hi {q} \times 
    \real \hi {q} \to \real \)
be positive definite kernels, and let $\ca H \lo U \hii 0 $ and $\ca H \lo V \hii 0$ be the RKHS's generated by $\ka \lo U$ and $\ka \lo V$, respectively. As discussed in \cite{li2017} and \cite{li2018sufficient}, it suffices to focus on the subspaces of $\ca H \lo U \hii 0 $ and $\ca H \lo V \hii 0$ spanned by $\{ \ka \lo U (\cdot, u) - E [\ka \lo U ( \cdot, U)]: u \in \Omega \lo U \}$ and $\{ \ka \lo V (\cdot, v) - E [\ka \lo V ( \cdot, V)]: v \in \Omega \lo V \}$, respectively, where, for example, $E [\ka \lo U ( \cdot, U)]$ denotes the Bochner integral $\int \ka \lo U (\cdot, u) \, d P \lo U (u)$. We refer to these subspaces as $\ca H \lo U$ and $\ca H \lo V$.  \cite{li2017} and \cite{li2018sufficient} showed that the orthogonal complements of $\ca H \lo U$ and $\ca H \lo V$ are spanned by constants, which do not affect conditional independence. 
Let $\ca H \lo U \otimes \ca H \lo V$ denote the completion of the inner product space
\begin{align}\label{eq:span}
    \spn \{ f g: f \in \ca H \lo U, g \in  \ca H \lo V \}, 
\end{align}
where $f g$ represents the function $(u,v) \mapsto f(u) g(v)$, with its inner product determined by 
\begin{align}\label{eq:inn prod}
    \langle f \lo 1 g \lo 1, f \lo 2 g \lo 2 \rangle \lo {\ca H \lo U \otimes \ca H \lo V} = \langle f \lo 1, f \lo 2 \rangle \lo {\ca H \lo U} \, \langle g \lo 1, g \lo 2 \rangle \lo {\ca H \lo V}. 
\end{align}
An alternative (but equivalent) definition of the tensor product space $\ca H \lo U \otimes \ca H \lo V$ is the Hilbert space generated by the tensor products $f \otimes g $, which are linear operators from $\ca H \lo V$ to $\ca H \lo U$ defined by $(f \otimes g) (h) = f \langle g, h \rangle \lo {\ca H \lo V} $. This formulation was also adopted in \cite{tang2024nonparametric}, but for the current development, the definition of $\ca H \lo U \otimes \ca H \lo V$ via (\ref{eq:span}) and (\ref{eq:inn prod}) is more convenient. 

For each $x \in \Omega \lo X$, it admits the decomposition \( x = \sum \lo {i=1} \hi r U \lo i (x) V \lo i \trans (x)\), where 
\begin{align*}
\ali     U \lo 1 (x) = \lambda \lo 1 (x) \hi {1/2} U \lo 1 \hi 0 (x) \;, \ldots, \; U \lo r (x) = \lambda \lo r (x) \hi {1/2} U \lo r \hi 0 (x), \\
\ali     V \lo 1 (x) = \lambda \lo 1 (x) \hi {1/2} V \lo 1 \hi 0 (x) \;, \ldots, \; V \lo r (x) = \lambda \lo r (x) \hi {1/2} V \lo r \hi 0 (x). 
\end{align*}
We then define the feature map of $X$ as \( F: \Omega \lo X \to \ca H \lo U \otimes \ca H \lo V, \), where 
\begin{align}\label{eq:feature map}
\begin{split}
F (X) = \ali   \sum \lo {i=1} \hi r [ \ka \lo U (\cdot, U \lo i (X)) - E  \ka \lo U (\cdot, U \lo i (X)) ] \, [\ka \lo V (\cdot , V \lo i(X)) - E \ka \lo V ( \cdot, V \lo i (X))]. 
\end{split}
\end{align}
A similar  feature was used by \cite{virta2024} in their tensor   preserving PCA.

\subsection{Tucker- and CP-form nonlinear SDR for tensors}

In this subsection, we develop two approaches to nonlinear SDR for tensors, in parallel with their linear counterparts.
\vspace{.1in}

\noindent
{\em Tucker form. } We assume that there exist $f \lo 1, \ldots, f \lo s \in \ca H \lo U$ and $g \lo 1, \ldots, g \lo t \in \ca H \lo V$ such that 
\begin{align}\label{eq:nonlinear tucker}
     Y\independent X|\{\langle f_k g \lo \ell , F (X)   \rangle_{\ca H \lo U \otimes \ca H \lo V } : k =1,\cdots,s, \, \ell=1,\cdots,t \}.   
\end{align}
Henceforth, we abbreviate $\langle \cdot, \cdot \rangle \lo {\ca H \lo U \otimes \ca H \lo V}$ as $\langle  \cdot, \cdot \rangle \lo \otimes$.  
Under this notation, we can rewrite  $\langle f \lo k g \lo \ell , F (X) \rangle \lo \otimes$  as 
\begin{align*}
\ali \sum \lo {i=1} \hi r     \langle f \lo k, \ka \lo U (\cdot, U \lo i (X)) \rangle \lo {\ca H \lo U} \, \langle g \lo \ell, \ka \lo V (\cdot , V \lo i(X)) \rangle \lo {\ca H \lo V }= \sum \lo {i=1} \hi r      f \lo k (   U \lo i (X))  \,  g \lo \ell(  V \lo i(X)). 
\end{align*}
When $f_k(U_i(X))$ and $g_\ell(V_i(X))$ are linear functions, that is, $f_k(U_i(X)) = \alpha_k^\top U_i(X)$ and $g_\ell(V_i(X)) = \beta_\ell^\top V_i(X)$, the right-hand side above reduces to
$\sum_{i=1}^r \alpha_k^\top U_i(X) V_i(X)^\top \beta_\ell = \alpha_k^\top X \beta_\ell$. Hence, (\ref{eq:nonlinear tucker}) reduces to the sufficient dimension folding problem in \cite{Li2010}. 
For convenience, we refer to this problem as the Tucker-form nonlinear tensor SDR (Tucker-NTSDR).

\medskip

\noindent{\em CP form.} In this case, we assume that there exist orthonormal sets $\{f_1, \ldots, f_d\}$ in $\ca H \lo U$ and $\{g_1, \ldots, g_d\}$ in $\ca H \lo V$ such that
\begin{align}\label{eq:nonlinear cp}
    Y \indep X | \{ \langle f \lo 1 \,  g \lo 1, F (X)\rangle \lo {\otimes}, \ldots, \langle f \lo d \, g \lo d, F(X)  \rangle \lo \otimes \}. 
\end{align}
Again, when $f \lo k$ and $g \lo k$ are linear functions, the above dimension reduction problem reduces to (\ref{eq:linear cp}), the linear  SDR for tensor introduced by \cite{zhong2015}.  We refer to problem (\ref{eq:nonlinear cp}) as the CP-form nonlinear tensor SDR (CP-NTSDR).

\subsection{Tensor envelopes} 
As in the linear case, there is also a ``stretched-out'' version of tensor nonlinear SDR. 
Let $A_1, \ldots, A_d$ be elements of $\ca H \lo U \otimes \ca H \lo V$. 
The ``stretched-out'' version of tensor nonlinear SDR is  defined through the conditional independence
\begin{align}\label{eq:nonlinear stretched}
    Y \independent X | \langle A \lo 1 , F (X) \rangle \lo \otimes, \ldots,   \langle A \lo d, F (X) \rangle \lo \otimes.  
\end{align}
Let $\| \cdot \| \lo \otimes$ denote the norm induced by the inner product $\langle \cdot, \cdot \rangle \lo \otimes$. 
Each $A_i$ in the above relation is the $\| \cdot \| \lo \otimes$-limit of a sequence whose elements are of the form 
\( c_1 f_1 g_1 + \cdots + c_k f_k g_k, \)
where $c_1, \ldots, c_k \in \real$, $f_1, \ldots, f_k \in \ca H \lo U$, and $g_1, \ldots, g_k \in \ca H \lo V$. 
When $f_1, \ldots, f_k$ and $g_1, \ldots, g_k$ are linear functions, the function 
\(\langle c_1 f_1 g_1 + \cdots + c_k f_k g_k, F(X) \rangle \lo {\otimes}\)
reduces to
\begin{align*}
(c \lo 1 (\beta \lo 1 \otimes \alpha \lo 1) + \cdots + c \lo k ( \beta \lo k \otimes \alpha \lo k)) \trans \vec (X). 
\end{align*}
Since vectors of the form $c \lo 1 (\beta \lo 1 \otimes \alpha \lo 1) + \cdots + c ( \beta \lo k \otimes \alpha \lo k)$ span the entire $\real \hi {pq}$ space, problem (\ref{eq:nonlinear stretched}) reduces to (\ref{eq:linear stretch}) in the linear case. 
We refer to problem (\ref{eq:nonlinear stretched}) as the stretched-out NTSDR. 

If either the Tucker-form NTSDR (\ref{eq:nonlinear tucker}) or the CP-form NTSDR (\ref{eq:nonlinear cp}) is satisfied, then it always corresponds to a stretched-out NTSDR, but the sufficient predictors may not be the most parsimonious. 
Specifically, for the Tucker-form NTSDR, we can write
\begin{align*}
\{\langle f \lo k \, g \lo \ell , F (X)\rangle \lo  \otimes : k = 1, \ldots, s, \ell = 1, \ldots, t \} = \{ \langle \tilde A \lo u, F (X) \rangle \lo \otimes: u = 1, \ldots, st  \}, 
\end{align*}
where $\tilde A \lo 1 = f \lo 1 g \lo 1$, $\tilde A \lo 2 = f \lo 1 g \lo 2$, \ldots, $\tilde A \lo {st} = f \lo s g \lo t$. 
Thus (\ref{eq:nonlinear tucker}) is equivalent to 
\begin{align*}
Y \independent X | \{ \langle   \bar A \lo u,  F(X)  \rangle \lo \otimes: u = 1, \ldots, st \}. 
\end{align*}
It will be shown that, if both (\ref{eq:nonlinear tucker}) and (\ref{eq:nonlinear stretched}) are   satisfied, then 
\begin{align*}
   \spn \{   A \lo 1, \ldots, A \lo d \} \subseteq \spn \{ f \lo k \otimes g \lo l : \, k = 1, \ldots, s, l = 1, \ldots, t\}. 
\end{align*}
Similarly, for the CP-form NTSDR, we can rewrite
(\ref{eq:nonlinear cp}) as
\begin{align*}
Y \independent X | \{  \langle f \lo k \otimes g \lo k,  F(X)  \rangle \lo \otimes: k=1, \ldots, c  \}. 
\end{align*}
It will be shown that, if both (\ref{eq:nonlinear cp}) and (\ref{eq:nonlinear stretched}) are both satisfied, then 
\begin{align*}
    \spn \{ A \lo 1, \ldots, A \lo d \} \subseteq \spn \{ f \lo k \otimes g \lo k : \, k = 1, \ldots, c \}. 
\end{align*}
In both cases, the goal is to identify the smallest subspaces in $\ca H \lo U \otimes \ca H \lo V$ of the Tucker form or the CP form that contain the stretched-out NTSDR space spanned by $\{A_1, \ldots, A_d\}$. 
We refer to these spaces as the Tucker envelope and the CP envelope, respectively. In the following sections, we develop these three spaces in detail and their implementation methods.

\section{Tensor Regression Operator}
In this section, we extend GSIR to the stretched-out NTSDR setting. We introduce the covariance and regression operators in this context and show that, under mild conditions, the range of the regression operator is contained within the stretched-out NTSDR central subspace.

\subsection{Covariance operators involving the tensor feature $F(X)$} 

We now introduce the mean elements and covariance operators in the tensor product space \( \mathcal{H}_U \otimes \mathcal{H}_V \). 
The most convenient way to define these objects is through the Bochner integral; see, for example, \cite{hsing2015theoretical}. 
Let $(\Omega, \mathcal{F}, P)$ be a probability space and let $\mathcal{B}$ be a Banach space. 
Denote by $\mathcal{F}(\mathcal{B})$ the Borel $\sigma$-field generated by the open sets of $\mathcal{B}$. 
Let $f: \Omega \to \mathcal{B}$ be a function that is measurable with respect to $\mathcal{F} / \mathcal{F}(\mathcal{B})$. 
Since $\| f \|_{\mathcal{B}}$ is a nonnegative measurable function from $\Omega$ to $\mathbb{R}$, its integral 
\( E \| f \|_{\mathcal{B}} = \int_{\Omega} \| f \|_{\mathcal{B}} \, dP \)
is well defined. 
If this integral is finite, then the Bochner integral of $f$, denoted by \( \int_{\Omega} f \, dP \), is a well-defined  nonrandom element in $\mathcal{B}$. 
We will use the Bochner integral to define the aforementioned mean elements and covariance operators.

First, for a random matrix $X: \Omega \to \Omega_X$, the feature map $F(X)$ is the measurable function
\[
    F(X): \Omega \to \mathcal{H}_U \otimes \mathcal{H}_V, \quad \omega \mapsto F(X(\omega)).
\]
Hence, $F(X)$ is a random element of $\mathcal{H}_U \otimes \mathcal{H}_V$. 
Under the condition that \( E \| F(X) \|_{\otimes} < \infty \), \( E[F(X)] \) is a well-defined  nonrandom element of \( \mathcal{H}_U \otimes \mathcal{H}_V \). Next, consider the space 
\begin{align}\label{eq:tensor tensor}
    (\ca H \lo U \otimes \ca H \lo V) \otimes (\ca H \lo U \otimes \ca H \lo V ). 
\end{align}
This tensor product space, unlike the tensor product space $\ca H \lo U \otimes \ca H \lo V$ itself, is defined according to the following rule. Let $\ca H \lo 1$ and $\ca H \lo 2$ be two generic Hilbert spaces. Let $h \lo 1 \in \ca H \lo 1$ and $h \lo 2 \in \ca H \lo 2$. The tensor product $h \lo 1 \otimes h \lo 2$ is defined as the linear operator from $\ca H \lo 2 $ to $\ca H \lo 1$ such that, for any $g \in \ca H \lo 2$, $(h \lo 1 \otimes h \lo 2) g = h \lo 1 \langle h  \lo 2, g \rangle \lo {\ca H \lo 2}$. Consider the vector space 
\begin{align*}
 \ca A =    \spn \{ h \lo 1 \otimes h \lo 2: \, h \lo 1 \in \ca H \lo 1, h \lo 2 \in \ca H \lo 2 \}. 
\end{align*}
For any $h \lo 1 \otimes h \lo 2, g \lo 1 \otimes g \lo 2 \in \ca H \lo 1 \otimes \ca H \lo 2$, define the inner product between them as 
\begin{align}\label{eq:tensor inner product}
    \langle h \lo 1 \otimes h \lo 2, g \lo 1 \otimes g \lo 2 \rangle \lo {\ca H \lo 1 \otimes \ca H \lo 2} = \langle h \lo 1, g \lo 1 \rangle \lo {\ca H \lo 1} 
\, \langle h \lo 2, g \lo 2 \rangle \lo {\ca H \lo 2}.
\end{align}
This rule uniquely determines an inner product on $\mathcal{A}$ by the bilinearity of the inner product. 
Equipped with this inner product, $\mathcal{A}$ becomes an inner product space. 
The tensor product space $\mathcal{H}_1 \otimes \mathcal{H}_2$ is then defined as the completion of this inner product space. 
Note that this definition of the tensor product is --- at least in appearance --- different from the one we used for $\mathcal{H}_U \otimes \mathcal{H}_V$. 
Nevertheless, $\mathcal{H}_U \otimes \mathcal{H}_V$ could have been defined in the same way, and there exists a mapping that makes the two definitions equivalent. 
For convenience in the subsequent development, we define the tensor product space $(\mathcal{H}_U \otimes \mathcal{H}_V) \otimes (\mathcal{H}_U \otimes \mathcal{H}_V)$ according to the above rule, while retaining the earlier definition of $\mathcal{H}_U \otimes \mathcal{H}_V$.

Having defined (\ref{eq:tensor tensor}), and with $F(X)$ being a random element in $\mathcal{H}_U \otimes \mathcal{H}_V$, we see that $F(X) \otimes F(X)$ is a random element in (\ref{eq:tensor tensor}). 
Its Bochner integral with respect to $P$,
\[
    E[F(X) \otimes F(X)] = \int_{\Omega} F(X) \otimes F(X) \, dP,
\]
is well defined if
\begin{align}\label{eq:EFF}
E \left[ \| F(X) \otimes F(X) \|_{(\mathcal{H}_U \otimes \mathcal{H}_V) \otimes (\mathcal{H}_U \otimes \mathcal{H}_V)} \right] < \infty,
\end{align}
where the expectation on the left is the usual expectation for random numbers. By the tensor inner product rule (\ref{eq:tensor inner product}), \(    E \| F (X) \otimes F(X) \| \lo {(\ca H \lo U \otimes \ca H \lo V) \otimes ( \ca H \lo U \otimes \ca H \lo V )} =     E ( \| F (X)  \| \lo { \ca H \lo U \otimes \ca H \lo V }\hi 2 ). \)

Hence, under the condition \( E \| F(X) \|_{\mathcal{H}_U \otimes \mathcal{H}_V}^2 < \infty \), 
the Bochner expectation in (\ref{eq:EFF}) is well defined. 
We now define the linear operator $\Sigma_{FF}$ as
\begin{align}\label{equ:Covariance}
    \Sigma_{FF} = E[F(X) \otimes F(X)] - E[F(X)] \otimes E[F(X)].
\end{align}

Note that $\Sigma_{FF}$ is a nonrandom element of $(\mathcal{H}_U \otimes \mathcal{H}_V) \otimes (\mathcal{H}_U \otimes \mathcal{H}_V)$. 
Let $\kappa_Y : \Omega_Y \times \Omega_Y \to \mathbb{R}$ be a positive definite kernel, and let $\mathcal{H}_Y^{0}$ denote the RKHS generated by $\kappa_Y$. 
Following the constructions of $\mathcal{H}_U$ and $\mathcal{H}_V$, we define $\mathcal{H}_Y$ as the subspace of $\mathcal{H}_Y^{0}$ spanned by
\[
\{ \kappa_Y(\cdot, y) - E[\kappa_Y(\cdot, Y)] : \, y \in \Omega_Y \}.
\]

Let $(\mathcal{H}_U \otimes \mathcal{H}_V) \otimes \mathcal{H}_Y$ be the tensor product space of $\mathcal{H}_U \otimes \mathcal{H}_V$ and $\mathcal{H}_Y$, defined in the same manner as in (\ref{eq:tensor tensor}); that is, an element of $(\mathcal{H}_U \otimes \mathcal{H}_V) \otimes \mathcal{H}_Y$ is a linear operator from $\mathcal{H}_Y$ to $\mathcal{H}_U \otimes \mathcal{H}_V$. 
The Bochner integrals \( E [\ka \lo Y (\cdot, Y)] \mbox{ and } E [F(X) \otimes (\ka \lo Y (\cdot, Y) - E \ka \lo Y (\cdot, Y) ) ]\) are well defined if 
$E\| \kappa_Y(\cdot, Y) \|_{\mathcal{H}_Y} < \infty $
and $E\| F(X) \otimes \kappa_Y(\cdot, Y) \|_{(\mathcal{H}_U \otimes \mathcal{H}_V) \otimes \mathcal{H}_Y} < \infty$.

 We denote the linear operator $E [F(X) \otimes (\ka \lo Y (\cdot, Y) - E \ka \lo Y (\cdot, Y) ) ]$
by  $\Sigma \lo {FY}$. The next assumption summarizes the conditions under which $\Sigma \lo {FF}$ and $\Sigma \lo {FY}$ are well defined.

\begin{assumption}\label{assumption:integrable features}
\quad 
\( E \| F(X) \|_{\mathcal{H}_U \otimes \mathcal{H}_V}^2 < \infty \)
and
\( E[\kappa_Y(Y, Y)] < \infty. \)
\end{assumption}

\noindent Both conditions are mild and hold for a wide class of kernels. The following proposition provides sufficient conditions for the first part of Assumption~\ref{assumption:integrable features}.

\begin{proposition}\label{ka2ka2}
Suppose that, for each $i = 1, \ldots, r$, 
\[E \left[ \kappa_U(U_i(X), U_i(X))^2 \right] < \infty \text{ and } E \left[ \kappa_V(V_i(X), V_i(X))^2 \right] < \infty.\]
Then \( E \| F(X) \|_{\mathcal{H}_U \otimes \mathcal{H}_V}^2 < \infty. \)
\end{proposition}

\subsection{Regression operator on tensor feature} As in non-tensorial nonlinear SDR frameworks \citep{li2011psvm, li2017, li2018sufficient}, 
the regression operator also plays a central role in tensor nonlinear SDR. 
Before defining this operator, we first introduce several notations and concepts regarding the kernel, range, domain, and the Moore--Penrose generalized inverse of a linear operator. 

For two Hilbert spaces $\ca G$ and $\ca H$ and a linear operator $A: \ca G \to \ca H$, the kernel of $A$, denoted by $\ker (A)$,  is the set $\{ x \in \ca G: A x = 0 \}$. The range of $A$, denoted by $\ran (A)$, is the set $\{ A x: x \in \ca G \}$. The domain of $A$ is denoted by $\dom (A)$. While $\ker (A)$ is always a closed linear subspace,   $\ran (A)$ is a linear subspace that may not be closed. We use $\cran (A)$ to denote the closure of $\ran(A)$. If $\ca H = \ca G$ and  $A$ is a self-adjoint operator, then $\ker (A) = \cran(A) \oc$. For a subset $\ca S$ of $\ca H$, we use $A|\ca S$ to denote the restriction of $A$ on $\ca S$; that is, $A|\ca S$ is a mapping from $\ca S \to \ca H$ such that, for any $h \in \ca S$, $(A | \ca S)(h) = A (h)$. 
Suppose  $A: \ca H \to \ca H$  is a  self-adjoint operator. Since $\cran (A) = \ker(A) \oc$,   $\ker(A | \cran (A))=\{0\}$, which  implies that the mapping $(A|\cran(A)): \cran (A) \to \ran(A)$ is an injection. We define $(A|\cran(A))\inv: \ran(A) \to \cran(A)$ as the \emph{Moore-Penrose generalized inverse} of $A$, and denote it by $A \hi \dagger$.

Note that if $\ran (\Sigma \lo {FY}) \subseteq \ran ( \Sigma \lo {FF})$, then \(\ran(\Sigma_{FY}) \subseteq \dom(\Sigma_{FF}^\dagger)\), 
and
\(\Sigma_{FF}^\dagger \Sigma_{FY}\) 
is a well-defined linear map from 
\(\mathcal{H}_Y\) 
to 
\(\cran(\Sigma_{FF})\). 
We refer to this operator as the \emph{tensor regression operator}. The formal assumption and definition are given as follows. 

\begin{assumption}\label{assumption:ran in ran} \quad 
    $\ran (\Sigma \lo {FY}) \subseteq \ran ( \Sigma \lo {FF})$. 
\end{assumption}

\begin{definition} \label{definition:regression}
Under Assumptions \ref{ass:assumptions2},  \ref{assumption:integrable features}, and \ref{assumption:ran in ran}, the linear operator $\Sigma \lo {FF} \hi \dagger \Sigma \lo {FY}$ is called the \emph{tensor regression operator}, and is denoted by $R \lo {FY}$.     
\end{definition}

\section{Generalized Slice Inverse Regression for stretched-out NTSDR} The critical step in the subsequent development is to establish a connection between the tensor regression operator \(R_{FY}\) and the central $\sigma$-field. 
Let \(\mathcal{H}_X + \mathbb{R}\) denote the set \(\{\, h + c : h \in \mathcal{H}_X,\, c \in \mathbb{R} \,\}\). We begin by introducing two notations for conditional expectation. 
Let \(S\) and \(T\) be two random elements defined on the probability space \((\Omega, \mathcal{F}, P)\), taking values in measurable spaces \((\Omega_S, \mathcal{F}_S)\) and \((\Omega_T, \mathcal{F}_T)\), respectively. One may view the conditional expectation $E(S|T)$ as a mapping from \(\Omega\) to \(\Omega_{S}\) that is measurable with respect to \(\sigma(T)\). 
Since $E(S|T)$ is \(\sigma(T)\)-measurable, there exists a measurable function \(h : \Omega_T \to \Omega_S\) such that \(E(S|T) = h \circ T\). 
Alternatively, one may interpret \(E(S|T)\) as the function \(h\) itself. 
Although both interpretations appear in the literature, it is seldomly necessary to distinguish between them. However, our context requires  this distinction for clarity. We will denote the first version as $E(S|T)$ and the second as $E(S|T = \cdot)$. Thus, $E(S|T)$ is a mapping from $\Omega$ to $\Omega \lo {S}$; whereas $E(S|T = \cdot)$ is a mapping from $\Omega \lo T$ to $\Omega\lo S$. We use $P \lo X$, $P \lo Y$, and $P \lo {XY}$ to denote the distributions of $X$, $Y$, and $(X,Y)$, respectively. That is, $P \lo X = P \of X \inv$, $P \lo Y = P \of Y \inv$, and $P \lo {XY} = P \of (X,Y) \inv$. Let 
\begin{align*}
\ca H \lo X =   \{  x \mapsto  \langle f, F(x) \rangle \lo \otimes: \, f \in \ca H \lo U \otimes \ca H \lo V \}. 
\end{align*}
That is, $\ca H \lo X$ is the collection of all functions of $x$ of the form $x \mapsto \langle f, F(x) \rangle \lo \otimes$ where $f \in \ca H \lo U \otimes \ca H \lo V$.  Clearly, $\ca H \lo X$ is a linear space. For any $\tilde f, \tilde g \in \ca H \lo X$ defined by $\tilde f(x) = \langle f, F(x) \rangle \lo \otimes$ and $\tilde g(x) = \langle g, F(x) \rangle \lo \otimes$, we define their inner product as 
$ \langle \tilde f, \tilde g \rangle \lo {\ca H \lo X} = \langle f, g \rangle \lo {\otimes}. $
It follows immediately that $\ca H \lo X$ and $\ca H \lo U \otimes \ca H \lo V$ are isomorphic.

\begin{theorem}\label{theorem:fisher consistency 1}
Suppose Assumptions \ref{ass:assumptions2}--\ref{assumption:ran in ran} and the following conditions are satisfied:
    \begin{enumerate}
        \item $\ka \lo Y: \Omega \lo Y \times \Omega \lo Y \to \real $ is a universal kernel; 
        \item for each $g \in \ca H \lo Y$, $E[g(Y)|X = \cdot \, ] \in \ca H \lo X + \real$;  
        \item $R \lo {FY}$ is a finite-rank operator. 
    \end{enumerate}
Then $\sigma \{ h(X): h \in \ran (R \lo {FY}) \}$ is the central $\sigma$-field. 
\end{theorem}

As will be seen in the proof of this theorem in the Supplementary Material, the third condition that the tensor regression has finite rank is not crucial, as it can be removed with careful handling of the closure of the range of $R \lo {FY}$. However, for meaningful dimension reduction, it is entirely reasonable to make this assumption, and this simplifies the proof.

\section{Tucker Tensor Envelope and its Estimation} \label{section:tucker form envelope}
 
\subsection{Tucker envelope}

\def\lat{\mathrm{Lat}}
\def\latuv{\lat (\ca H \lo U) \otimes \lat (\ca H \lo V)}

Let $\lat (\ca H \lo U)$ and $\lat (\ca H \lo V)$ denote the collections of all subspaces of $\ca H \lo U$ and $\ca H \lo V$, respectively. 
Here, the symbol $\lat$ stands for ``lattice,'' as these collections of subspaces indeed form lattices (see, for example, \cite{conway2019functional}, p.~178). 
Define 
\begin{align*}
\lat (\ca H \lo U) \otimes \lat (\ca H \lo V) 
= \{ \ca S \lo U \otimes \ca S \lo V : \ca S \lo U \in \lat (\ca H \lo U), \, \ca S \lo V \in \lat (\ca H \lo V) \}. 
\end{align*}
The following proposition shows that the intersection of two members of $\latuv$ remains within $\latuv$. 
Its proof is straightforward and thus omitted. 

\begin{proposition}\label{proposition:intersection}
If $\ca S \lo U \otimes \ca S \lo V$ and $\ca S \lo U' \otimes \ca S \lo V'$ are members of $\latuv$, then $(\ca S \lo U \otimes \ca S \lo V) \cap (\ca S \lo U' \otimes \ca S \lo V') 
= (\ca S \lo U \cap \ca S \lo U') \otimes (\ca S \lo V \cap \ca S \lo V'). $   
\end{proposition}

This property motivates the following definition of the Tucker envelope.

\begin{definition}
The intersection of all members of $\lat ( \ca H \lo U) \otimes \lat ( \ca H \lo V) $ that contain $\ran (R \lo {FY})$ is the \emph{Tucker envelope} of $\ran ( R \lo {FY} )$, which we denote by $\ca E \lo T (\ran (R \lo {FY}))$. 
\end{definition}

\subsection{Population-level objective function for Tucker envelope} To guide the sample-level implementation, we next introduce an optimization criterion at the population level whose solution yields the Tucker envelope of $\ran (R \lo {FY})$. 
We make the following assumption. 

\begin{assumption}\label{assumption:finitedimension}
The Tucker envelope is a finite-dimensional subspace.  
\end{assumption}

\def\bigtimes{\mbox{\Large{$\times$}}}

\noindent 
Suppose $\ca E \lo T ( \ran (R \lo {FY})) = \ca S \lo U \otimes \ca S \lo V$. 
We refer to $\ca S \lo U$ as the \emph{left Tucker NTSDR space} and to $\ca S \lo V$ as the \emph{right Tucker NTSDR space}. 
Let the dimensions of $\ca S \lo U$ and $\ca S \lo V$ be $s$ and $t$, respectively. 
The following theorem shows that, at the population level, the Tucker envelope can be obtained by solving a sequence of least-squares problems. In what follows, define
\begin{align*}
\ca H \lo U \hi s =  \overset{s}{\underset{i=1}{\bigtimes}} \ca H \lo U, \quad
\ca H \lo V \hi t =  \overset{t}{\underset{j=1}{\bigtimes}} \ca H \lo V, \quad     
\ca H \lo Y \hi {s \times t} =  \overset{s}{\underset{i=1}{\bigtimes}} \overset{t}{\underset{j=1}{\bigtimes}} \ca H \lo Y, \quad 
\ca H \lo {YUV} =   \ca H \lo Y \hi {s \times t} \times \ca H \lo U \hi s \times \ca H \lo V \hi t,  
\end{align*}
where $\bigtimes$ denotes the Cartesian product. 
Let $f = (f \lo 1, \ldots, f \lo s)$, $g = (g \lo 1, \ldots, g \lo t)$, and $h = \{ h \lo {ij}: i = 1, \ldots, s,\, j = 1, \ldots, t \}$.

\def\leftnorm{\left\|\phantom{\int}\right.\hspace{-.15in}} 
\def\rightnorm{\hspace{-0.15in}\left.\phantom{\int}\right\|}

\begin{theorem}\label{thm:theorem7}
Let $(f \hi *, g \hi *, h \hi *) \in \ca H \lo {YUV}$ be the solution to the optimization problem:
\begin{align}\label{equ:equation1}
    \min_{(f,g,h) \in \ca H \lo {YUV}} 
    E \leftnorm R \lo {FY} \,  \ka \lo Y 
    - \sum_{i=1}^{s}\sum_{j=1}^{t} f_i  g_j h_{ij}(Y) 
    \rightnorm_{\ca H \lo U \otimes \ca H \lo V}^2.
\end{align}
Then,
\begin{align*}
\spn\{ f \lo 1 \hi *, \ldots, f \lo s \hi * \} 
\otimes 
\spn \{ g \lo 1 \hi *, \ldots, g \lo t \hi * \}  
= \ca E \lo T ( \ran (R \lo {FY})).
\end{align*}
\end{theorem}

\def\leftnorm{\left\|\phantom{\int}\right.\hspace{-.15in}} 
\def\rightnorm{\hspace{-0.15in}\left.\phantom{\int}\right\|}
\subsection{Explicit solution at population level} The objective function in (\ref{equ:equation1}) is essentially a least-squares criterion, 
which can be minimized iteratively, with each step admitting an explicit solution. 
The iterative procedure consists of the following three steps:
\begin{enumerate}
    \item[\textbf{Step 1.}] Given $f \in \ca H \lo U \hi s$ and $g \in \ca H \lo V \hi t$, 
    minimize (\ref{equ:equation1}) with respect to $h \in \ca H \lo Y \hi {s \times t}$;
    \item[\textbf{Step 2.}] Given $g \in \ca H \lo V \hi t$ and $h \in \ca H \lo Y \hi {s \times t}$, 
    minimize (\ref{equ:equation1}) with respect to $f \in \ca H \lo U \hi s$;
    \item[\textbf{Step 3.}] Given $h \in \ca H \lo Y \hi {s \times t}$ and $f \in \ca H \lo U \hi s$, 
    minimize (\ref{equ:equation1}) with respect to $g \in \ca H \lo V \hi t$. 
    \item[]  
    Repeat Steps~1--3 until convergence.
\end{enumerate}

The next two theorems and a corollary provide explicit solutions for each of these three steps. 
To derive these solutions, we first introduce a linear operator mapping from $\ca H \lo V$ to $\ca H \lo U$, 
associated with the random element $R \lo {FY}[ \ka \lo Y (\cdot, Y) - E \ka \lo Y (\cdot, Y)]$ in $\ca H \lo U \otimes \ca H \lo V$. 
For a given element $r \in \ca H \lo U \otimes \ca H \lo V$, consider the bilinear form \(B: \, \ca H \lo U \times \ca H \lo V \to \real, B(f,g) = \langle f g, r \rangle \lo {\otimes}\) where $f \in \ca H \lo U$ and $g \in \ca H \lo V$. By the Cauchy–Schwarz inequality and the inner product rule (\ref{eq:tensor inner product}),
\begin{align*}
    B(f,g) \le \| fg \| \lo \otimes \, \| r \| \lo \otimes = \| r \| \lo \otimes \, \| f \| \lo {\ca H \lo U} \, \| g \| \lo {\ca H \lo V}. 
\end{align*}
This shows that $B$ is bounded bilinear form. 
By Theorem~2.2 of \cite{conway2019functional}, there exists a unique bounded linear operator 
$T \lo r: \ca H \lo V \to \ca H \lo U$ such that, for any 
$f \in \ca H \lo U$ and $g \in \ca H \lo V$, 
\begin{align}\label{eq:T r}
    \langle f, T \lo r (g) \rangle \lo {\ca H \lo U} 
    = \langle T \lo r \hi * (f), g \rangle \lo {\ca H \lo V} 
    = \langle f g, r \rangle \lo \otimes,  
\end{align}
where $T \lo r \hi *$ is the adjoint operator of $T \lo r$. 
In this way, every element of $\ca H \lo U \otimes \ca H \lo V$ can be associated with a unique linear operator 
from $\ca H \lo V$ to $\ca H \lo U$. 
A particular element of interest is 
$s(Y) = R \lo {FY} [\ka \lo Y (\cdot, Y) -  E  \ka \lo Y (\cdot, Y)${, for which the} operator 
$T \lo {s(Y)}: \ca H \lo V \to \ca H \lo U$ 
plays a key role in deriving the explicit solutions for the iterative optimization.

\def\tsum{\textstyle{\sum}}

Let $\ca H$ be a generic Hilbert space, and let   
$m = (m \lo 1, \ldots, m \lo k) \trans$ and 
$n = (n \lo 1, \ldots, n \lo \ell )\trans$ 
be vectors whose components are elements of $\ca H$. 
We use  $\langle m, n \trans \rangle \lo {\ca H}$ 
to denote the $k \times \ell$ matrix whose $(i,j)$th entry is 
$\langle m \lo i, n \lo j \rangle \lo {\ca H}$. We will also use the matrix product between a matrix of real numbers and a matrix (or vector) whose entries are elements of a Hilbert space. 
Specifically, let $A \in \real \hi {a \times b}$ be a matrix with real-valued entries, 
and let $g = (g \lo 1, \ldots, g \lo b) \trans$ be a vector whose components belong to a Hilbert space $\ca H$. 
At first glance, the product $A g$ may appear to be undefined since the entries of $A$ and $g$ lie in different spaces. 
However, the notation is legitimate if we interpret it in terms of scalar multiplication and vector addition in $\ca H$: 
we define 
$A g = \sum \lo {i=1} \hi a \sum \lo {j=1} \hi b a \lo {ab} g \lo b$, 
where $a \lo {ab}$ is the $(a,b)$th entry of $A$. 
This sum is well defined because, under the operations of scalar multiplication and addition in a vector space, 
it represents a valid element of $\ca H$. Finally, we use $H(y)$ to denote the $s \times t$ matrix $\{ h \lo {ij}(y) \}$, and $h(y)$ to denote its vectorization, $\vec (H(y))$.

\begin{theorem}\label{theorem:explicit sol for step1}  For given $h \in \ca H \lo Y \hi {s \times t}$ and $g \in \ca H \lo V \hi t$, the minimizer of (\ref{equ:equation1}) over $f \in \ca H \lo U \hi s$ is 
\begin{align}\label{eq:f star}
    f \hi *  = \{ E  [H(Y) \langle g, g \trans   \rangle \lo {\ca H \lo V} H(Y) \trans ]  \} \inv E  [ H(Y)  T \lo {s(Y)}   g    ]. 
\end{align}
\end{theorem}

Using the same reasoning, one can obtain the explicit solution for Step~2. 
The result is summarized in the following corollary, whose proof is omitted for brevity.  

\begin{corollary}\label{corollary:explicit sol for step2}  
For given $h \in \ca H \lo Y \hi {s \times t}$ and $f \in \ca H \lo U \hi s$, 
the minimizer of (\ref{equ:equation1}) over $g \in \ca H \lo V \hi t$ is 
\begin{align*}
    g \hi *  
    = \{ E [H (Y) \trans \langle f,  f \trans \rangle \lo {\ca H \lo U} H(Y) ] \}\inv 
      E  [ H(Y) \trans  T \lo {s(Y)} \hi * f  ].
\end{align*}
\end{corollary}

Next, we minimize the objective function in (\ref{equ:equation1}) with respect to 
$h_{ij}(y)$, for fixed $f \in \ca H \lo U \hi s$ and $g \in \ca H \lo V \hi t$. 
In what follows, the notation $\otimes \lo k$ indicates the Kronecker product between matrices. 
Specifically, for an $m \times n$ matrix $A = \{ a \lo {ij} \}$ and an arbitrary matrix $B$, we define
\begin{align*}
    A \otimes \lo k B = 
    \begin{pmatrix}
        a \lo {11} B & \cdots & a \lo {1n} B \\
        \vdots & \ddots & \vdots \\
        a \lo {m1} B & \cdots & a \lo {mn} B
    \end{pmatrix}.
\end{align*}
Note that the Kronecker product $\otimes \lo k$ is different from the tensor product $\otimes$. 
We allow the factors of $\otimes \lo k$ to be matrices whose entries are elements of a Hilbert space. 
We denote the minimizer of (\ref{equ:equation1}) with respect to $h$ by $h \hi *$.

\begin{theorem}\label{theorem:explicit sol for step3}  
For given $f \in \ca H \lo U \hi s$ and $g \in \ca H \lo V \hi t$, 
the minimizer of the objective function in (\ref{equ:equation1}) over 
$h \in L \lo 2 (P \lo Y)  \hi {s \times t}$ is 
\begin{align}\label{eq:h star}
   h \hi * (y)  
   = 
   ( \langle g, g \trans \rangle \lo {\ca H \lo V} \inv 
   \otimes \lo k 
   \langle f, f \trans \rangle \lo {\ca H \lo U} \inv )   
   \,  \langle g  \otimes \lo k f , s(y) \rangle \lo \otimes.
\end{align}
\end{theorem}

\section{CP Tensor Envelope}\label{section:cp form envelope} 
Unlike the Tucker envelope, which regresses $s(Y)$ on a function of the form 
$\sum \lo {i=1} \hi s \sum \lo {j=1} \hi t h \lo {ij} f \lo i g \lo j$, 
the CP envelope successively regresses $s(Y)$ on functions of the   rank-one form 
$u(y) f g$, where $f \in \ca H \lo U$ and $g \in \ca H \lo V$, 
until no variation remains unexplained. 
Specifically, we consider the objective function:
$E (     \| s (Y) -  u(Y) \, f g\| \lo {\otimes} \hi 2 )$, where $f \in \ca H \lo U$, $g \in \ca H \lo V$, and $u  \in L  \lo 2 (P \lo Y)$. Let 
\begin{align*}
    L(f,g) = \min \lo {u \in L \lo 2 (P \lo Y)} \, E (     \| s (Y) -  u(Y) \, f g\| \lo {\otimes} \hi 2 ). 
\end{align*}

We now introduce the CP-sequence of $s(Y)$ and the CP-envelope of $\ran (R \lo {FY})$. 
Let $\ca E^{C}(s(Y)) = \spn\{ f_1 g_1, \ldots, f_d g_d \}$ 
denote the CP envelope of $s(Y)$ obtained through the minimization process 
described in the following definition.

\begin{definition}\label{definition:cp envelope}\em 
The \emph{CP-sequence} of $s(Y)$ is the sequence $\{(f \lo k, g \lo k): k = 1, 2, \ldots \}$ in $\ca H \lo U \times \ca H \lo V$, 
defined recursively through the following iterative minimization:
\begin{enumerate}
    \item Step~1: minimize $L(f,g)$ to obtain $(f \lo 1, g \lo 1)$;
    \item Step~$k$ ($k = 2, 3, \ldots$): minimize $L(f,g)$ over
    \begin{align*}
        (f, g) \in (\spn \{f \lo 1, \ldots, f \lo {k-1}\}) \oc 
        \times 
        (\spn \{g \lo 1, \ldots, g \lo {k-1}\}) \oc,
    \end{align*}
    where $(f \lo 1, g \lo 1), \ldots, (f \lo {k-1}, g \lo {k-1})$ are the solutions obtained in Steps~1 through $k-1$.
\end{enumerate}
\end{definition}

\noindent 
Intuitively, the CP-sequence extracts one pair of directions $(f_k, g_k)$ at a time, 
each capturing the dominant remaining component of $s(Y)$ in $\ca H \lo U \otimes \ca H \lo V$. This motivates the following assumption, which guarantees that the sequence terminates in finitely many steps.

\begin{assumption}\label{assumption:CP5}
In the minimization process of Definition~\ref{definition:cp envelope}, 
there exists a finite integer $k \ge 1$ such that 
$L(f \lo \ell, g \lo \ell) = 0$ for all $\ell = k, k + 1, \ldots$.
\end{assumption}

\noindent 
The next theorem establishes that this space is indeed the smallest rank-one–structured subspace 
that contains the range of the regression operator $R \lo {FY}$.

\def\qh{\red{(question here) }}

\begin{theorem}\label{theorem:cp contains} 
Suppose that Assumption~\ref{assumption:CP5} holds and  
\begin{align*}
\Sigma  \lo {YY} = E[(\ka \lo Y (\cdot, Y) - E [\ka \lo Y (\cdot, Y)])\otimes (\ka \lo Y (\cdot, Y) - E [\ka \lo Y (\cdot, Y)])]
\end{align*}
is a nonsingular operator. 
If
\begin{align*}
    d = \min \{ k: L(f \lo \ell, g \lo \ell) = 0 \ \text{for all } \ell = k, k + 1, \ldots \},  
\end{align*}
then 
$\spn \{f \lo 1 g \lo 1, \ldots, f \lo d g \lo d \}$ 
is the smallest linear space of the form 
$\spn \{f \lo 1 g \lo 1, \ldots, f \lo k g \lo k \}$ 
that contains $\ran (R \lo {FY})$. 
\end{theorem}

This theorem motivates the following definition of the CP-envelope. 

\begin{definition}
Under the conditions of Theorem~\ref{theorem:cp contains}, 
the linear subspace $\spn \{ f \lo 1 g \lo 1, \ldots, f \lo d g \lo d \}$ is called the \emph{CP envelope} of $\ran (R \lo {FY})$.
\end{definition}

\noindent
Theorem~\ref{theorem:cp contains} further implies that 
$\spn \{ f \lo 1 g \lo 1, \ldots, f \lo d g \lo d \}$ 
is the smallest linear subspace of the form 
$\spn \{ f \lo 1 g \lo 1, \ldots, f \lo k g \lo k \}$ 
that satisfies
\begin{align*}
    Y \indep X \,\big|\,
    \{ \langle f \lo 1 g \lo 1, F(X) \rangle \lo \otimes, 
       \ldots, 
       \langle f \lo d g \lo d, F(X) \rangle \lo \otimes \}.
\end{align*}
Like the Tucker envelope, each step of the minimization in Definition~\ref{definition:cp envelope} 
can be implemented through an iterative sequence of three least-squares problems, 
each of which admits an explicit analytical solution. 
The following corollary is a special case of 
Theorem~\ref{theorem:explicit sol for step1}, 
Corollary~\ref{corollary:explicit sol for step2}, 
and Theorem~\ref{theorem:explicit sol for step3}, 
with $s = t = 1$. 
Its proof is omitted, as it follows directly from the preceding results. 
Define
\begin{align*}
    L \lo C (f, g, u) 
    = E \big( \| s(Y) - u(Y) f g \| \lo \otimes \hi 2 \big),
\end{align*}
where the subscript $C$ in $L \lo C$ indicates the ``CP-approach.''

\begin{corollary}\label{corollary:explicit sol for cp}  \ \  
\begin{enumerate}
   \item For given $u \in L \lo 2 (P \lo Y)$ and $g \in \ca H \lo V$, the minimizer of $L \lo C(f,g,u)$ over $f \in \ca H \lo U $ is 
\begin{align*}
    f \hi *  = E  [u (Y)  T \lo {s(Y)}   g    ] /    E  [u \hi 2 (Y) \| g \| \hi 2 \lo {\ca H \lo V} ]. 
\end{align*}
\item  For given $u \in L \lo 2 (P \lo Y) $ and $f \in \ca H \lo U$, the minimizer of $L \lo C(f, g, u)$ over $g \in \ca H \lo V$ is 
\begin{align*}
    g \hi *  = E  [u (Y)  T \lo {s(Y)} \hi *   f    ] /    E  [u \hi 2 (Y) \| f \| \hi 2 \lo {\ca H \lo U} ]. 
\end{align*}
\item For given $f \in \ca H \lo U $  and $g \in \ca H \lo V $, the minimizer of $L \lo C (f,g,u)$ over $u \in L \lo 2 (P \lo Y)$  is 
\begin{align*} 
   u \hi * (y)  = 
    \langle fg , s(y) \rangle \lo \otimes / \| f g \| \lo \otimes \hi 2 . 
\end{align*}
\end{enumerate}
\end{corollary}

Alternatively, as the next theorem shows,  the above minimization problem can be re-expressed as what resembles an eigenvalue problem, which can, in fact, be solved by a sequence of iterative eigenvalue problems. Let 
\begin{align*}
    M(f, g) = E (  \langle f, T \lo {s(Y)} \, g \rangle \lo {\ca H \lo U} \hi 2 ). 
\end{align*}

\begin{proposition}\label{proposition:resemble eigen}
The minimization problem in Definition~\ref{definition:cp envelope} 
is equivalent to the following sequence of maximization problems:
\begin{enumerate}
    \item Step~1: maximize $M(f,g)$ subject to 
    $\| f \|_{\ca H \lo U} = 1$ and $\| g \|_{\ca H \lo V} = 1$; 
    \item Step~$k = 2, \ldots, d$: maximize $M(f,g)$ subject to 
    \begin{align*}
        \| f \|_{\ca H \lo U} = 1, \quad 
        \langle f, f_i \rangle_{\ca H \lo U} = 0, \ i = 1, \ldots, k-1, \\[0.4em]
        \| g \|_{\ca H \lo V} = 1, \quad 
        \langle g, g_i \rangle_{\ca H \lo V} = 0, \ i = 1, \ldots, k-1,
    \end{align*}
    where $f_1, \ldots, f_{k-1}$ are the solutions obtained in 
    Steps~1 through $k-1$.
\end{enumerate}
\end{proposition}

\noindent
Since, for any bounded linear operator $A: \ca H \lo V \to \ca H \lo U$, 
\(\langle f, A g \rangle_{\ca H_U}^2
= \langle f, A (g \otimes g) A^* f \rangle_{\ca H_U} \), we can re-express $M(f,g)$ as
\[
M(f,g)
= \langle f, E [T_{s(Y)} (g \otimes g) T_{s(Y)}^*] f \rangle_{\ca H_U}
= \langle g, E [T_{s(Y)}^* (f \otimes f) T_{s(Y)}] g \rangle_{\ca H_V}.
\]
Therefore, for a fixed $g$, the optimal $f^*$ satisfies the eigen-equation
\[
E[(T_{s(Y)} g) \otimes (T_{s(Y)} g)] f^* = \alpha f^*,
\]
and symmetrically, for a fixed $f$, the optimal $g^*$ satisfies the eigen-equation
\[
E[(T_{s(Y)}^* f) \otimes (T_{s(Y)}^* f)] g^* = \beta g^*.
\]
These eigen-equations are equivalent to the least-squares solutions 
in Corollary~\ref{corollary:explicit sol for cp}, 
with
\[
\alpha = E[\langle f^* g, s(Y) \rangle_{\otimes}^2],
\qquad
\beta  = E[\langle f g^*, s(Y) \rangle_{\otimes}^2].
\]

\section{Asymptotic Theory} In this section, we establish the consistency and convergence rates of the proposed nonlinear tensor SDR estimator, 
focusing on the CP-version of NTSDR. 
The proof of the following lemma follows the same line of reasoning as Lemma~5 in \cite{fukumizu2007statistical} 
and is therefore omitted for brevity.

\begin{lemma}\label{lem:bound}
Under Assumptions~\ref{ass:assumptions2}, \ref{assumption:integrable features}, and \ref{assumption:ran in ran}, 
the covariance operators \( \Sigma_{FF} \) and \( \Sigma_{FY} \) defined in Section~3.1 are Hilbert--Schmidt. Moreover, 
\[
\|\hat{\Sigma}_{FY} - \Sigma_{FY}\|_{\mathrm{HS}} = O_P(n^{-1/2}), 
\qquad 
\|\hat{\Sigma}_{FF} - \Sigma_{FF}\|_{\mathrm{HS}} = O_P(n^{-1/2}).
\]
\end{lemma}

To obtain the convergence rate of \( R_{FY} \), we impose the following additional smoothness condition linking the covariance operators.

\begin{assumption}\label{assumption:smoothness}
There exists a bounded linear operator \( S_{FY} \) and a constant \( 0 < \beta \le 1 \) such that
$\Sigma_{FY} = \Sigma_{FF}^{1+\beta} S_{FY}$.
\end{assumption}

This assumption is similar to that used in  \cite{li2017}, 
and it ensures that the image of \( \Sigma_{FY} \) is sufficiently concentrated in the low-frequency region 
of the spectrum of \( \Sigma_{FF} \), thereby representing a degree of smoothness in the relation between \( X \) and \( Y \); 
see \cite{li2017} and \cite{li2018sufficient} for further discussions on this point. 
Define
\[
\hat{R}_{FY} = (\hat{\Sigma}_{FF} + \epsilon_n I)^{-1} \hat{\Sigma}_{FY},
\]
where \( I: \ca H_U \otimes \ca H_V \to \ca H_U \otimes \ca H_V \) denotes the identity operator, 
and \( \epsilon_n \) is a tuning parameter. In the next theorem, we assume that \( S_{FY} \) is either a bounded linear operator or a Hilbert--Schmidt operator, 
corresponding to different convergence modes.

\begin{theorem}\label{thm:theorem17}
Suppose  Assumptions~\ref{ass:assumptions2}, \ref{assumption:integrable features}, \ref{assumption:ran in ran}, and \ref{assumption:smoothness} are satisfied. 
    \begin{enumerate}
        \item If \( S_{FY} \) is a bounded operator, then \ 
       $ \| \hat R \lo {FY} - R \lo {FY} \|_{\mathrm{OP}} = O_P(\epsilon \lo n  \hi \beta + \epsilon \lo n \inv n^{-1/2}).$
        \item If \( S_{FY} \) is a Hilbert-Schmidt operator, then \ 
        $ \| \hat R \lo {FY} - R \lo {FY} \|_{\mathrm{HS}} = O_P(\epsilon \lo n  \hi \beta + \epsilon \lo n \inv n^{-1/2}).$
    \end{enumerate}
    \end{theorem}

The proof follows the same reasoning as that of Theorem~5 in \cite{li2017} and is therefore omitted. 
Intuitively, it is natural to conjecture that the sufficient predictors achieve the same convergence rate 
of order \( \epsilon_n^{\beta} + \epsilon_n^{-1} n^{-1/2} \). 
However, establishing this rate requires substantial technical development. 
To avoid too much digression, in this paper we only establish the consistency of the sufficient predictors 
and leave the precise rate analysis to future research. 

Let \( s(Y) = R_{FY}[ \kappa_Y(\cdot, Y)- E \ka \lo Y (\cdot, Y)] \) and 
\( \hat{s}(Y) = \hat{R}_{FY} [\kappa_Y(\cdot, Y)- E \lo n \ka \lo Y (\cdot, Y)]  \) , 
and define
\[
R(t) = E\|s(Y) - t(Y)\|_{\otimes}^2, 
\quad 
R_n(t) = E_n\|s(Y) - t(Y)\|_{\otimes}^2, 
\quad 
\hat{R}_n(t) = E_n\|\hat{s}(Y) - t(Y)\|_{\otimes}^2.
\]
Consider the class of functions of the form
\begin{align*}
 \{  \tsum \lo {k=1} \hi d u \lo k(Y) f \lo k g \lo k: \ali   u \lo 1, \ldots, u \lo d \in L \lo 2 (P \lo Y),  \{f \lo 1, \ldots, f \lo d\} \mbox{ is an orthonormal set in $\ca H \lo U$}, \\
\ali \{g \lo 1, \ldots, g \lo d \} \mbox{ is an orthonormal set in $\ca H \lo V$} \}.
\end{align*} 
Let \( \ca F \) be a subset of this class, 
and define \( \ca U = \{ u(y) = \|s(y) - t(y)\|_{\otimes}^2 : t \in \ca F \} \). 
Let
\(\hat{t}(Y) = \tsum_{k=1}^d \hat{u}_k(Y) \hat{f}_k \hat{g}_k\)
be the minimizer of 
\( \hat{R}_n(t) \)
over \( \ca F \). 
Let $N \lo {[\,]}(\epsilon, \ca G, L \lo 1 (P) )$
denote the \( \epsilon \)-bracketing number of a class of functions \( \ca G \) 
(see, e.g., \cite{vaart1996weak}).

\begin{theorem}
\label{thm:consistency}
Suppose that the following conditions hold:
\begin{enumerate}
    \item $\{(X_i, Y_i)\}_{i=1}^n$ are i.i.d.\ observations from the joint distribution of $(X, Y)$;
    \item for every $\epsilon > 0$, the $\epsilon$-bracketing number satisfies 
    \( N_{[\, ]}(\epsilon, \mathcal{G}, L \lo 1(P)) < \infty \);
    \item $E\{\kappa_Y(Y, Y)\} < \infty$;
    \item $s(Y) \in \mathcal{F}$.
\end{enumerate}
Then the estimated sufficient predictor is consistent in the sense that
\[
\|\hat{t}(Y) - s(Y)\|_{\mathcal{F}} \stackrel{P}{\longrightarrow} 0.
\]
\end{theorem}

\section{Sample-Level Implementation}
\subsection{Coordinate Mapping}

The optimization procedures in Sections~\ref{section:tucker form envelope} and~\ref{section:cp form envelope} are implemented at the sample level by representing elements and linear operators in Hilbert spaces through coordinate mappings.  In this subsection, we outline the basic definitions and rules of coordinate representation (or coordinate mapping). Let $\ca H$ be a finite-dimensional Hilbert space spanned by $\ca B=\{b_1,\ldots,b_m\}$. Any $f\in\ca H$ can be written as $f=\sum_{k=1}^m c_k b_k$, whose coordinate vector with respect to $\ca B$ is $[f]_{\ca B}=(c_1,\ldots,c_m)^\top$. Since
\(\langle f,b_k\rangle_{\ca H}=\sum_{j=1}^m c_j\langle b_j,b_k\rangle_{\ca H},\)
we obtain
\begin{align}\label{eq:explicit coordinate for function}
[f]_{\ca B}=G_{\ca B}^\dagger
\begin{pmatrix}
\langle f,b_1\rangle_{\ca H}\\
\vdots\\
\langle f,b_m\rangle_{\ca H}
\end{pmatrix},
\end{align}
where $G_{\ca B}=\{\langle b_i,b_j\rangle_{\ca H}\}_{i,j=1}^m$ is the Gram matrix of $\ca B$.

Let $\ca H_1$ and $\ca H_2$ be Hilbert spaces spanned by $\ca B_i=\{b_1^{(i)},\ldots,b_{n_i}^{(i)}\}$, $i=1,2$. For a linear operator $A:\ca H_1\to\ca H_2$, its coordinate matrix with respect to $(\ca B_1,\ca B_2)$ is
\(_{\ca B_2}[A]_{\ca B_1}
=([A b_1^{(1)}]_{\ca B_2},\ldots,[A b_{n_1}^{(1)}]_{\ca B_2}),\)
or equivalently,
\begin{align}\label{eq:coord op}
_{\ca B_2}[A]_{\ca B_1}
=G_{\ca B_2}^\dagger
\begin{pmatrix}
\langle A b_1^{(1)},b_1^{(2)}\rangle_{\ca H_2} & \cdots & \langle A b_{n_1}^{(1)},b_1^{(2)}\rangle_{\ca H_2}\\
\vdots & \ddots & \vdots\\
\langle A b_1^{(1)},b_{n_2}^{(2)}\rangle_{\ca H_2} & \cdots & \langle A b_{n_1}^{(1)},b_{n_2}^{(2)}\rangle_{\ca H_2}
\end{pmatrix}.
\end{align}
The mappings $f\mapsto[f]_{\ca B}$ and $A\mapsto\,_{\ca B_2}[A]_{\ca B_1}$ are called the \emph{coordinate mappings} of functions and linear operators. Following \cite{li2018nonparametric}, we use $\ca B(\ca H_1,\ca H_2)$ to denote the class of all linear operators from $\ca H_1$ to $\ca H_2$. The sets $\ca B \lo 1$ and $\ca B \lo 2$ are spanning system that need not be linearly independent. For a finite subset $\{h_1,\ldots,h_m\}\subset\ca H$, its Gram matrix is $G=\{\langle h_i,h_j\rangle_{\ca H}\}_{i,j=1}^m$. For any integer $n$, we denote $\{1, \ldots, n \}$ by  $[n]$. The following theorem, which will be used heavily in our derivations, is taken from \cite{li2018nonparametric}.

\begin{theorem} \label{thm:coord}
    Let $(\mathcal{H},\mathcal{B})$ and $\{(\mathcal{H}_i,\mathcal{B}_i):i=1,2,3\}$ be finite-dimensional Hilbert spaces with spanning systems $\mathcal{B}$, $\mathcal{B}_1$, $\mathcal{B}_2$, and $\mathcal{B}_3$ of sizes $d$, $d_1$, $d_2$, and $d_3$, respectively.  
Then the coordinate mappings satisfy the following properties: \\
    1. (Evaluation) For any $h\in \mathcal{H}_1$ and $ A  \in \mathcal{B}(\mathcal{H}_1,\mathcal{H}_2)$, $[ A h]_{\mathcal{B}_2} = (_{\mathcal{B}_2}[ A ]_{\mathcal{B}_1})[h]_{\mathcal{B}_1}$.\\
    2. (Linearity) If $ A _1, A _2 \in \mathcal{B}(\mathcal{H}_1,\mathcal{H}_2)$ and $c_1,c_2 \in \mathbb{R}$, then  
    \begin{align*}
      _{\mathcal{B}_2}[c_1 A _1+c_2 A _2]_{\mathcal{B}_1} = c_1(_{\mathcal{B}_2}[ A _1]_{\mathcal{B}_1})+c_2(_{\mathcal{B}_2}[ A _2]_{\mathcal{B}_1}). 
    \end{align*}
    3. (Composition) If $ A _1\in \mathcal{B}(\mathcal{H}_1,\mathcal{H}_2)$ and $ A _2 \in \mathcal{B}(\mathcal{H}_2,\mathcal{H}_3)$, then
\begin{align*}
  _{\mathcal{B}_3}[ A _2 A _1]_{\mathcal{B}_1} = (_{\mathcal{B}_3}[ A _2]_{\mathcal{B}_2})(_{\mathcal{B}_2}[ A _1]_{\mathcal{B}_1}).  
\end{align*}
       4. (Inner Product) If $h_1,h_2 \in \mathcal{H} $, then $\langle h_1,h_2 \rangle_{\mathcal{H} } = ([h_1]_{\mathcal{B} })\trans G_{\mathcal{B}} ([h_2]_{\mathcal{B}})$, where $G_{\mathcal{B}}$ is the Gram matrix of $\ca B$. \\
    5. (Identity) If $I \in \mathcal{B}(\mathcal{H} ,\mathcal{H} )$ is the identity operator, then $_{\mathcal{B} }[I]_{\mathcal{B} } = Q_{\mathcal{B} }$, where $Q_{\mathcal{B} }$ is the projection onto $\operatorname{span}\{[b_1]_{\mathcal{B} },\cdots,[b_{d}]_{\mathcal{B} }\}$ in $\ca H$.\\
    6. (Tensor Product) If $h_1\in \mathcal{H}_1$ and $h_2 \in \mathcal{H}_2$, then $_{\mathcal{B}_2}[h_2\otimes h_1]_{\mathcal{B}_1} = [h_2]_{\mathcal{B}_2}[h_1]_{\mathcal{B}_1}\trans G_{\mathcal{B}_1} $, where $G_{\mathcal{B}_1}$ is the Gram matrix for $\ca B \lo 1$.
\end{theorem}

\subsection{Selection of Bases and Construction of Gram Matrices} To construct $\ca H_U$, let
\begin{align*}
f_{ia}^0 = \ka_U(\cdot, U_i(X_a)), \quad  
f_{ia} = f_{ia}^0 - E_n(f_i^0),
\end{align*}
where $E_n(f_i^0) = n^{-1}\sum_{a=1}^n f_{ia}^0$.  
Then $\ca H_U$ is the linear span of $\{f_{ia}: i\in[r], a\in[n]\}$.  
Similarly, $\ca H_V$ is spanned by $\{g_{ia}: i\in[r], a\in[n]\}$ where
$g_{ia}^0=\ka_V(\cdot, V_i(X_a))$ and $g_{ia}=g_{ia}^0-E_n[g_i^0]$, and  
 $\ca H_Y$ is spanned by $\{h \lo a: a \in [n] \}$, where $h_a=h_a^0-E_n [h^0]$ and $h_a^0=\ka_Y(\cdot,Y_a)\}$.

The tensor-product space $\ca H_U\otimes\ca H_V$ is spanned by
\(\ca B_{UV}=\{f_{ia}g_{jb}:i,j\in[r],a,b\in[n]\},\)
which has $(rn)^2$ elements.  
However, $s(Y_a)$ lies in a much smaller subspace. By definition,
\[
F(X_a)-E_nF(X)=\sum_{i=1}^r(r_{ia}-E_n r_i), \quad 
r_{ia}=f_{ia}g_{ia},\quad
E_n r_i=n^{-1}\sum_{a=1}^n f_{ia}g_{ia},
\]
so $F(X_a)$ lies in the subspace spanned by $\{r_{ia}:i\in[r],a\in[n]\},$ 
which has only $rn$ members.  
Hence, the coordinate representation of $F(X_a)$ can be constructed within this $rn$-dimensional subspace—a key computational gain of the tensor-respecting map $F(X)$.

\def\one{\mathbbm{1}}

In summary, the spanning systems for the four spaces are
\begin{align*}
\ca B_U &= \{f_{ia}:i\in[r],a\in[n]\}, \quad
\ca B_V = \{g_{ia}:i\in[r],a\in[n]\},\\
\ca B_{UV} &= \{f_{ia}g_{jb}:i,j\in[r],a,b\in[n]\}, \quad
\ca B_Y = \{h_a:a\in[n]\}.
\end{align*}

We next compute the corresponding Gram matrices.  
For any integer $m$, let $Q_m=I_m-{\mathbbm{1}}_m{\mathbbm{1}}_m^\top/m$,  
where $I_m$ is the $m\times m$ identity and ${\mathbbm{1}}_m$ is the $m$-vector of ones.  
Note that $Q_m$ is the projection onto the orthogonal complement of $\spn({\mathbbm{1}}_m)$ in $\real^m$.  
The following lemma summarizes a useful property of Gram matrices in tensor-product spaces.

\begin{lemma}\label{lemma:gram tensor}
    Suppose that $\ca H$ is a Hilbert space spanned by $\ca B = \{ b \lo 1, \ldots, b \lo m \}$ and $\ca G$ is a Hilbert space spanned by $\ca C = \{ c \lo 1, \ldots, c \lo n \}$. Let $\ca D = \{ b \lo i c \lo j: i \in [m], j \in [n] \}$, and assume $b \lo i c \lo j$ is ordered so that $j$ runs the first and $i$ runs the second: that is, the members of $\ca D$ are ordered as $b \lo 1 c \lo 1, \ldots, b \lo 1  c \lo n, \cdots \cdots,     b \lo m c \lo 1, \ldots, b \lo m  c \lo n$. Let $G \lo {\ca B}, G \lo {\ca C}$, and $G \lo {\ca D}$ be the Gram matrices of the spanning systems $\ca B$, $\ca C$, and $\ca D$. Then \begin{align}\label{eq:G lo D}
G \lo {\ca D}=G \lo {\ca B} \otimes \lo k G \lo {\ca C}. 
\end{align}
\end{lemma}
Lemma~\ref{lemma:gram tensor} provides the tensor-product structure of Gram matrices, 
which directly facilitates the derivation of the following result for the specific bases constructed above.

\begin{theorem}\label{cor:gram}
The Gram matrices for   \(\ca B \lo U, \ca B \lo V, \ca B \lo {UV}\) and \(\ca B \lo Y\) are,  respectively, 
\begin{align*}
    G \lo U  = \ali (I \lo r \otimes \lo k Q \lo n) K \lo U (I \lo r \otimes \lo k  Q \lo n) ,\quad     G \lo V = (I \lo r \otimes \lo k  Q \lo n) K \lo V (I \lo r \otimes  \lo k 
    Q \lo n), \\
G \lo Y = \ali  Q \lo n K \lo Y Q \lo n, \quad     G_{UV} = G \lo U \otimes \lo k  G \lo V,
\end{align*}
where 
\begin{align*}
    K \lo U =\ali  \{ \ka \lo U ( U \lo i(X \lo a), U \lo j(X \lo b) ): (i,a) \in [r] \times [n], (j,b) \in [r] \times [n]\}, \\ 
    K \lo V = \ali  \{ \ka \lo V ( V \lo i(X \lo a), V \lo j(X \lo b) ): (i,a) \in [r] \times [n], (j,b) \in [r] \times [n]\}, \\
    K \lo Y =\ali \{\ka \lo Y (Y \lo a, Y \lo b): a,b \in [n] \}. 
\end{align*}
\end{theorem}

\subsection{Coordinate Representations of Linear Operators}

We now express the sample estimates of $\Sigma \lo {FF}$, $\Sigma \lo {FY}$, $R \lo {FY}$, $s(Y \lo a)$ in coordinate forms,
 relative to the spanning systems defined previously. These estimates are
\begin{align*}
\Sigma_{FF} &= E_n[F(X)\otimes F(X)] - E_n[F(X)]\otimes E_n[F(X)], \\
\Sigma_{FY} &= E_n[F(X)\otimes \kappa_Y(\cdot,Y)] - E_n[F(X)]\otimes E_n[\kappa_Y(\cdot,Y)], \\
R_{FY} &= \Sigma_{FF}^{-1}\Sigma_{FY}, \quad
s(Y \lo a ) = R_{FY}\big[\kappa_Y(\cdot,Y \lo a ) - E_n[\kappa_Y(\cdot,Y)]\big].
\end{align*}
For simplicity, we have used the same symbols $\Sigma_{FF}$, $\Sigma_{FY}$, and $R_{FY}$ for their sample-level estimates, omitting the hats. Let $e_{ia,jb}$ denote an $(nr)^2$-dimensional canonical basis vector whose $(i'a',j'b')$th entry equals $1$ if $(ia,jb)=(i'a',j'b')$ and $0$ otherwise.  
Throughout, the composite index $jb$ runs first and $ia$ second; within each $ia$, the index $a$ runs first ($a=1,\ldots,n$) and $i$ second ($i=1,\ldots,r$).  Let $\delta \lo a$ denote the $(nr) \hi 2$ dimensional vector $\sum \lo {i=1} \hi r e \lo {ia,ia}$, and let $\Delta$ denote the $(nr) \hi 2 \times n$ matrix $(\delta \lo 1, \ldots, \delta \lo n)$. These definitions will be used to construct coordinate expressions for $F(X_a)$ and $s(Y_a)$ in the following development.

\begin{proposition}\label{prop:sampcov}
The sample-level covariance and cross-covariance operators have the following coordinate representations:
\begin{enumerate}
    \item $\!_{\ca B_Y}[\Sigma_{YY}]_{\ca B_Y} = n^{-1} G_Y$;
    \item $\!_{\ca B_{UV}}[\Sigma_{FY}]_{\ca B_Y} = n^{-1} \Delta G_Y$;
    \item $\!_{\ca B_{UV}}[\Sigma_{FF}]_{\ca B_{UV}} = n^{-1} \Delta Q_n \Delta^\top G_{UV}$.
\end{enumerate}
\end{proposition}

To derive the coordinate form of $s(Y_a)$, we first establish two useful algebraic identities concerning $\Delta$ and $Q_n$ as defined above.

\begin{lemma}\label{eq:Delta T Delta}
$\Delta^\top \Delta = r I_n$ and $(\Delta Q_n \Delta^\top)^\dagger = r^{-2}\Delta Q_n \Delta^\top$.
\end{lemma}

This result is particularly useful because $\Delta Q_n \Delta^\top$ is an $(nr)^2\times(nr)^2$ matrix, which is extremely large in practice. 
Part~2 of Lemma~\ref{eq:Delta T Delta} implies that its Moore--Penrose inverse can be obtained analytically, 
thus avoiding explicit numerical inversion of this high-dimensional matrix. 
We next derive the coordinate representation of $s(Y_a)$ for any $a\in[n]$.

\begin{corollary}\label{cor:sampsY}
For each $a\in[n]$, the coordinate representation of $s(Y_a)$ with respect to $\mathcal{B}_{UV}$ is
$[s(Y_a)]_{\mathcal{B}_{UV}}
= r^{-1}(G_{UV})^\dagger \Delta G_Y e_a$.
\end{corollary}

We now develop the coordinate forms of the operator $T_{s(Y_a)}$ and its adjoint for any $a\in[n]$. 
For a vector $v=(v_k)$ with $k\in[mn]$, define $\mat_m(v)$ as the $m\times n$ matrix obtained by stacking the first $m$ entries of $v$ as the first column, the next $m$ entries as the second column, and so on. In other words, $\mat_m(\cdot)$ is the inverse of the $\vec(\cdot)$ operation.

\begin{corollary}\label{cor:sampIRFO}
The coordinate representation  of   \( T_{s(Y \lo a)} \) with respect to \(\mathcal{B}_U,\mathcal{B}_V\) is 
\begin{equation}\label{eq:coord IRFO}
    \lo {\mathcal{B}_U}  [T_{s(Y \lo a )}] \lo {\mathcal{B}_V} 
    =
r \hi {-1} G \lo {U} \hi \dagger    (G \lo U G \lo U \hi \dagger) \,  [\mat \lo {rn}     ( \Delta   G \lo Y e \lo a )]\trans \,  (G \lo V G \lo V \hi \dagger). 
\end{equation}
Similarly, the coordinate representation of the adjoint operator of $T \lo {s (Y \lo a)}$ is 
\begin{equation}\label{eq:coord IRFO adj}
    \lo {\mathcal{B}_V}  [T^{*}_{s(Y_a)}] \lo {\mathcal{B}_U} 
    =r \hi {-1}  G \lo {V} \hi \dagger    (G \lo V G \lo V \hi \dagger ) \,   \mat \lo {rn}    (\Delta   G \lo Y e \lo a ) \,  (G \lo U G \lo U \hi \dagger).     
\end{equation}
\end{corollary}

\def\mat{\mathrm{mat}}

It is well known that, for any symmetric matrix $A$, the product $A A^\dagger$ is the projection onto the column space of $A$. 
Since $G_U = (I_r \otimes_k Q_n) K_U (I_r \otimes_k Q_n)$ with $K_U$ nonsingular, the column space of $G_U$ coincides with that of $I_r \otimes_k Q_n$. 
Because $I_r \otimes_k Q_n$ itself is a projection matrix, we have 
\(G_U G_U^\dagger = I_r \otimes_k Q_n.\)
An identical relation holds for $G_V$. 
Consequently, the coordinate representations of the operators $T_{s(Y_a)}$ and $T^*_{s(Y_a)}$ admit the following forms.

\begin{corollary}\label{cor:alt sampIRFO}
The coordinate representation of the IRFO operator \(T_{s(Y_a)}\) is
\begin{equation}\label{IRFO}
\!_{\mathcal{B}_U}[T_{s(Y_a)}]_{\mathcal{B}_V}
= r^{-1} G_U^\dagger (I_r \otimes Q_n)
[\mat_{rn}(\Delta G_Y e_a)]^\top (I_r \otimes Q_n).
\end{equation}
Similarly, the coordinate representation of its adjoint operator is
\begin{equation}\label{IRFO_adj}
\!_{\mathcal{B}_V}[T^*_{s(Y_a)}]_{\mathcal{B}_U}
= r^{-1} G_V^\dagger (I_r \otimes Q_n)
\mat_{rn}(\Delta G_Y e_a) (I_r \otimes Q_n).
\end{equation}
\end{corollary}

\subsection{Sample-Level Estimation of Tucker Tensor Envelope}\label{subsection:sample tucker}

In this subsection, we implement the three-step iteration procedure in 
Theorem~\ref{theorem:explicit sol for step1}, 
Corollary~\ref{corollary:explicit sol for step2}, 
and Theorem~\ref{theorem:explicit sol for step3} at the sample level. 
Since vectors of Hilbert space-valued functions will be used repeatedly, 
we introduce the following notation.

\begin{definition} Suppose $\ca H$ is a Hilbert space spanned by $\ca B = \{b \lo 1, \ldots, b \lo m \}$ (again we allow members of $\ca B$ to be linearly dependent). Let $f = (f \lo 1, \ldots, f \lo s ) \trans \in \ca H \hi s$. Then $[f] \lo {\ca B}$ denotes the $s$ by $m$ matrix $([f \lo 1] \lo {\ca B}, \ldots, [f \lo s] \lo {\ca B} ) \trans$. 
\end{definition}

We now give the sample-level expressions of the three iterative updates for $(f,g,h)$.
The derivations of these expressions are given in Section S.2 of the Supplementary Material. Here we use $H(Y_a)$ to denote the $s\times t$ matrix 
$\{h_{ij}(Y_a)\}$.

\medskip
\noindent\textbf{Step 1: Updating $f$ given $(g,h)$.}  
Replacing the population expectations in~(\ref{eq:f star}) by the sample  averages, we have 
\begin{align*}
   [ f \hi * ] \lo {\ca B \lo U}  
= \ali  \left\{\tsum \lo {a=1} \hi n \left(  H (Y \lo a )   [g] \lo {\ca B \lo V} \, G \lo V \, [g] \lo {\ca B \lo V}\trans    H(Y \lo a )\trans  \right) \right\}\inv \\
\ali  r \hi {-1} \tsum \lo {a=1} \hi n  \left(  H(Y \lo a )    \left[  g \right] \lo {\ca B \lo V}   (G \lo V G \lo V \hi \dagger)   
 \,  \mat \lo {rn}    (\Delta   G \lo Y e \lo a )   \,   (G \lo U \hi \dagger G \lo U )    G \lo {U} \hi \dagger   \right). 
\end{align*}

\noindent\textbf{Step 2: Updating $g$ given $(f,h)$.}  
Similarly, the update for $g$ at the sample level is
\begin{align*}
   [ g \hi * ] \lo {\ca B \lo V}  
= \ali  \left\{\tsum \lo {a=1} \hi n \left(  H (Y \lo a ) \trans  [f] \lo {\ca B \lo U} \, G \lo U \, [f] \lo {\ca B \lo U}\trans    H(Y \lo a )  \right) \right\}\inv \\
\ali  r \hi {-1} \tsum \lo {a=1} \hi n  \left(  H(Y \lo a ) \trans   \left[  f \right] \lo {\ca B \lo U}   (G \lo U G \lo U \hi \dagger)   
 \,  \left( \mat \lo {rn}    (\Delta   G \lo Y e \lo a ) \right)\trans \,   (G \lo V \hi \dagger G \lo V )    G \lo {V} \hi \dagger   \right). 
\end{align*}

\noindent\textbf{Step 3: Updating $h$ given $(f,g)$.}  
Given $(f,g)$, the update for $h$ follows from Theorem~\ref{theorem:explicit sol for step3}:
\begin{align*}
h^*(Y_a)
= \ali 
r^{-1}
\Big\{
   ([g]_{\mathcal{B}_V} G_V [g]_{\mathcal{B}_V}^\top)^{\dagger}
   \otimes_k
   ([f]_{\mathcal{B}_U} G_U [f]_{\mathcal{B}_U}^\top)^{\dagger}
\Big\} \\
\ali 
\Big\{
  [ g\otimes_k f ]_{\mathcal{B}_{UV}}
 [ (G_U G \lo U \hi \dagger ) \otimes \lo k ( G_V G \lo V \hi \dagger ) ]
   \Delta G_Y e_a
\Big\},
\end{align*}
where 
\[
 [g\otimes_k f]_{\mathcal{B}_{UV}}
 =
\begin{pmatrix}
[f_1]_{\mathcal{B}_U}^\top\!\otimes_k\![g_1]_{\mathcal{B}_V}^\top\\
\vdots\\
[f_s]_{\mathcal{B}_U}^\top\!\otimes_k\![g_t]_{\mathcal{B}_V}^\top
\end{pmatrix}.
\]

\noindent\textbf{Tychonoff regularization.} For an $m \times m$ matrix positive semidefinite matrix $A$ and a positive number $\epsilon$, let $A(\epsilon)= A + \epsilon \lambda \lo {\max} (A) I \lo m $. 
In the above three steps, whenever we encounter $G_U^{\dagger}$ or $G_V^{\dagger}$,
we replace them by their Tychonoff-regularized forms $G \lo U (\epsilon \lo U)^{-1}$ and $G \lo V ( \epsilon \lo V)^{-1}$, 
where   $\epsilon_U,\epsilon_V>0$ are tuning constants (see Section~\ref{subsection:tuning}).
Similarly, the $G \lo U$ and $G \lo V$ in 
\begin{align*}
\ali \left( [g] \lo{\ca B \lo V} G \lo V  [g] \lo {\ca B \lo V} \trans \right) \hi   \dagger, \quad 
\left( [f] \lo{\ca B \lo U} G \lo U  [f] \lo {\ca B \lo U} \trans \right) \hi   \dagger,   \quad
  \left\{\tsum \lo {a=1} \hi n \left(  H (Y \lo a )   [g] \lo {\ca B \lo V} \, G \lo V \, [g] \lo {\ca B \lo V}\trans    H(Y \lo a )\trans  \right) \right\}\hi \dagger,  \\
\ali \hspace{1.4in} \left\{\tsum \lo {a=1} \hi n \left(  H (Y \lo a ) \trans  [f] \lo {\ca B \lo U} \, G \lo U \, [f] \lo {\ca B \lo U}\trans    H(Y \lo a )  \right) \right\}\hi \dagger, 
\end{align*}
are also replaced by $G \lo U (\epsilon \lo U)$ and $G \lo V (\epsilon \lo V)$. 
In addition, the coordinate representation of $s(Y_a)$ is regularized as
\begin{align}
\label{eq:coordinate of s(Y)}
[s(Y \lo b )]_{\mathcal{B}_{UV}} 
= \ali  r \hi {-1} [(G \lo U + \eta \lo U \lambda \lo {\max} (G \lo U) I \lo {rn} )  \inv  \otimes \lo k (G \lo V + \eta \lo V \lambda \lo {\max} (G \lo V) I \lo {rn} )  \inv  ]  \Delta    G_Y e \lo b \nonumber \\
\equiv  \ali  r \hi {-1} [ G \lo U \inv  (\eta \lo U)     \otimes \lo k  G \lo V \inv ( \eta \lo V )    ]  \Delta    G_Y e \lo b,
\end{align} 
where $\eta_U,\eta_V$ are separate tuning parameters. All these regularizations are applied consistently at each iteration 
to ensure numerical stability and prevent ill-conditioning 
in the sample-level updates of $(f,g,h)$.

\subsection{Sample-Level Estimation of CP Tensor Envelope}
\label{subsection:sample cp}

We now develop the estimation procedure for the CP tensor envelope. 
Let
\[
M(f,g)
= E\big[\langle f, T_{s(Y)}g\rangle_{\mathcal{H}_U}^2\big],
\]
which, at the sample level, is estimated by
\begin{align*}
    M(f, g) \equiv  \tsum \lo {a=1} \hi n \left(   [f] \lo {\ca B \lo U} \trans    \Gamma \lo a  [g] \lo {\ca B \lo V} \right) \hi 2, 
\end{align*}
where $\Gamma \lo a  =  n \inv r \hi {-2}      (G \lo U G \lo U \hi \dagger) \,  [\mat \lo {rn}     ( \Delta   G \lo Y e \lo a )]\trans \,  (G \lo V G \lo V \hi \dagger)$. 
The detailed derivation of the coordinate representation of
$\langle f, T_{s(Y_a)} g\rangle_{\mathcal{H}_U}$
and the resulting sample-level matrices is provided in the Supplementary Material Section S.3. By Proposition~\ref{proposition:resemble eigen}, the minimization of $M (f,g)$ can be turned into an iterative eigenvalue problem as follows. 

\noindent\textbf{Step 1. First pair of eigen vectors.}
Let $\Gamma_a$ be as defined above. 
For a given $[g]_{\mathcal{B}_V}$, compute $[f]_{\mathcal{B}_U}$ as the leading eigenvector of
\[
\sum_{a=1}^{n}\Gamma_a [g]_{\mathcal{B}_V}[g]_{\mathcal{B}_V}^\top \Gamma_a^\top.
\]
Given this $[f]_{\mathcal{B}_U}$, 
update $[g]_{\mathcal{B}_V}$ as the leading eigenvector of
\[
\sum_{a=1}^{n}\Gamma_a^\top [f]_{\mathcal{B}_U}[f]_{\mathcal{B}_U}^\top \Gamma_a.
\]
Iterate until convergence to obtain the first pair $(f_1^*, g_1^*)$.

\noindent\textbf{Steps 2 to \(d\). Subsequent pairs of eigenvectors.} 
After obtaining $(f_1^*,g_1^*),\ldots,(f_{k-1}^*,g_{k-1}^*)$, 
where the sets $\{f_i^*\}_{i=1}^{k-1}$ and $\{g_i^*\}_{i=1}^{k-1}$ are orthonormal in 
$\mathcal{H}_U$ and $\mathcal{H}_V$, respectively, 
define 
$P_{F_{k-1}}$ and $P_{G_{k-1}}$
as the projection matrices on to their spans; that is,
\begin{align*}
    F \lo {k-1} = ( [f \lo 1 \hi *] \lo {\ca B \lo U}, \ldots, [f \lo {k-1} \hi *] \lo {\ca B \lo U}), \quad 
    G \lo {k-1} = ( [g \lo 1 \hi *] \lo {\ca B \lo V}, \ldots, [ g \lo {k-1} \hi *] \lo {\ca B \lo V}), 
\end{align*}
Let
$Q_{F_{k-1}}=I_{rn}-P_{F_{k-1}}$ and 
$Q_{G_{k-1}}=I_{rn}-P_{G_{k-1}}$.
Update
\[
\Gamma_a \leftarrow Q_{F_{k-1}}^\top \Gamma_a Q_{G_{k-1}},
\]
and perform Step 1 on the updated $\Gamma \lo a$ and let $[\tilde f \lo k]\lo {\ca B \lo U}$ and $[\tilde g \lo k] \lo {\ca B \lo V}$ be the solutions. Then, 
\begin{align*}
[f \lo k \hi *]\lo {\ca B \lo U} = Q \lo {F \lo {k-1}}
[\tilde f \lo k]\lo {\ca B \lo U}, \quad 
[g \lo k \hi *] \lo {\ca B \lo V} = Q \lo {G \lo {k-1}}
[\tilde g \lo k] \lo {\ca B \lo V}. 
\end{align*}

\noindent\textbf{Step \(d+1\). Constructing the sufficient predictors.}
Once $\{(f_k^*,g_k^*)\}_{k=1}^d$ are obtained, we estimate
each sufficient predictor  as 
\begin{align*}
    \left([f \lo k \hi *] \lo {\ca B \lo U } \trans \otimes [g \lo k \hi *] \lo {\ca B \lo V} \trans \right) \ca B \lo {UV}, 
\end{align*}
where $\mathcal{B}_{UV}=\{f_{ia}g_{i'a'}\}$ denotes a vector 
with $a'$ running first, $i'$ second, $a$ third, and $i$ last.

The full coordinate representations and derivations 
for each step are detailed in the Supplementary Material Section S.3.
The iterative procedure has computational complexity $O(n^2 r^2)$,
making the CP decomposition an efficient and scalable alternative for tensor data.

\subsection{Selection of Tuning Parameters}\label{subsection:tuning}

\def\gcvt{\mathrm{GCV} \lo {\mathrm T}}
\def\gcvcp{\mathrm{GCV} \lo {\mathrm{CP}}} 

Proper selection of the tuning parameters ensures numerical stability and effective dimension reduction performance of the proposed method. 
These parameters include the kernel parameters in $\kappa_U$, $\kappa_V$, and $\kappa_Y$, 
as well as the Tychonoff regularization parameters $\epsilon_U$,  $\epsilon_V$, $\eta \lo U$, and $\eta \lo V$. For two pairs of indices $(i,a)$ and $(i',a')$, we write $(i,a) < (i',a')$ 
if either \(i < i'\) or \(i=i'\) and \(a < a'\). 
For the kernel bandwidth parameters $\gamma_U$, $\gamma_V$, and $\gamma_Y$, 
we adopt the following tuning procedure:
\[
\gamma_U = \frac{\rho_U}{2\sigma_U^2}, 
\qquad
\gamma_V = \frac{\rho_V}{2\sigma_V^2}, 
\qquad
\gamma_Y = \frac{\rho_Y}{2\sigma_Y^2},
\]
where the variance terms $\sigma_U^2$, $\sigma_V^2$, and $\sigma_Y^2$ are defined as:
\begin{align*}
\ali \sigma_U^2 = \binom{rn}{2}^{-1} \sum_{(i,a) < (i',a')}  \left\| U_{ia} - U_{i'a'} \right\|^2, \quad 
\sigma_V^2 = \binom{rn}{2}^{-1} \sum_{(i,a) < (i',a')}  \left\| V_{ia} - V_{i'a'} \right\|^2, \\
\ali \hspace{1.6in}
\sigma_Y^2 = {n \choose 2}^{-1}  \sum_{i<j}\|Y_i-Y_j\|^2.
\end{align*}
The parameters $\rho_U$, $\rho_V$, and $\rho_Y$ 
control the bandwidth of the corresponding kernels 
and are typically chosen within the range $\rho \in (0.1,10)$. Next, we develop a generalized cross validation procedure (GCV) to determine the Tychonoff regularization tuning parameters $\epsilon \lo U$ and $\epsilon \lo V$ for the Tucker and CP forms of NTSDR.
\begin{enumerate}
\item \textbf{Tuning parameters $(\eta_U, \eta_V)$ in $R_{FY}$}.
Let $\ka_Y^c(\cdot, y) = \ka_Y(\cdot, y) - E_n \ka_Y(\cdot, Y)$ denote the centered kernel. 
We select $(\eta_U, \eta_V)$ by minimizing a penalized least-squares criterion that matches 
$\langle s(Y_b), F(X_a) \rangle_{\otimes}$ with $\ka_Y^c(Y_a, Y_b)$ while preventing overfitting. 
The resulting generalized cross-validation (GCV) score takes the form
\begin{align}\label{eq:GCV-main}
\mathrm{GCV}(\eta_U, \eta_V)
= \frac{  \sum_{a,b=1}^n 
\big(\ka_Y^c(Y_a, Y_b) - \langle s(Y_b), F(X_a)\rangle_{\otimes}\big)^2}
{\big[\mathrm{tr}\{G_U(\eta_U)^{-1} G_U\}\,\mathrm{tr}\{G_V(\eta_V)^{-1} G_V\}-(nr)^2\big]^2}.
\end{align}
An explicit form of $\langle s(Y_b), F(X_a)\rangle_{\otimes}$ and the complete derivation of 
(\ref{eq:GCV-main}) are provided in Supplementary Material Section S.4.
\item \textbf{Tuning parameters $(\epsilon \lo U, \epsilon \lo V)$ for Tucker-form NTSDR.} \  Let  \(\{ (f \lo i \hi *, g \lo j \hi *, h \lo {ij} \hi * (Y)): i \in [s], j \in [t] \}\)
be the optimal solution described in Section \ref{subsection:sample tucker}. These solutions depend on the tuning parameters $\epsilon \lo U$ and $\epsilon \lo V$. 

Following the argument for the above construction, our GCV criterion in this case is
\begin{align*}
\mathrm{GCV} \lo {\mathrm{T}}(\epsilon \lo U, \epsilon \lo V) =     \frac{E \lo n \| s (Y) - \tsum \lo {i=1} \hi s \tsum \lo {j=1} \hi t f \lo i \hi * g \lo j \hi * h \lo {ij} \hi * (Y) \| \lo {\ca H \lo U \otimes \ca H \lo V} \hi 2 }{[ r \hi 2 (n-1) \hi 2 -\tr ( G \lo U G \lo U  \inv(\epsilon \lo U) ) \,  \tr ( G \lo V G \lo V \inv (\epsilon \lo V) )] \hi 2 }
\end{align*}

\item \textbf{Tuning parameters $(\epsilon \lo U, \epsilon \lo V)$    for CP-form NTSDR.} \  Let \(\{ (f \lo k \hi *, g \lo k \hi *, u \lo {k} \hi * (Y)): k \in [d] \} \) be the optimal solution described in Section \ref{section:cp form envelope}.   Our objective function for GCV in this case is 
\begin{align*}
\mathrm{GCV} \lo {\mathrm{CP}}(\epsilon \lo U, \epsilon \lo V) =     \frac{E \lo n \| s (Y) - \tsum \lo {k=1} \hi d f \lo k \hi * g \lo k \hi * u \lo k  \hi * (Y) \| \lo {\ca H \lo U \otimes \ca H \lo V} \hi 2 }{ [r \hi 2 (n-1) \hi 2 -\tr (   G \lo U G \lo U  \inv(\epsilon \lo U)  ) \,  \tr (  G \lo V G \lo V \inv (\epsilon \lo V)  )] \hi 2 }. 
\end{align*}
\end{enumerate}
We minimize  $\gcvt (\epsilon \lo U, \epsilon \lo V)$ or $\gcvcp (\epsilon \lo U, \epsilon \lo V)$ over a grid for $(\epsilon \lo U, \epsilon \lo V)$ to determine the optimal tuning constants.

\section{Nonlinear SDR for High-Order Tensors}
This section outlines an extension of our 
method to higher-order random tensors  \( X \in \mathbb{R}^{p_1 \times \cdots \times p_m} \), focusing on the CP-form NTSDR. The Tucker form NTSDR can be extended similarly, which is omitted.  

Our extension is based on  the CP-decomposition (or CANDECOMP/PARAFAC decomposition) of multiway tensors, which generalizes the singular value decomposition (SVD) for matrices. In the following, for vectors $a \lo 1, \ldots, a \lo r$, the tensor product $a \lo 1 \otimes \cdots \otimes a \lo r$ denotes the $r$-dimensional rank-1 tensor whose $(i \lo 1, \ldots, i \lo r)$th entry is $a \lo {i \lo 1} \cdots a \lo {i \lo r}$. Unlike the matrix SVD, which provides an exact factorization, the CP-decomposition for higher-order tensors is generally approximate. 
Specifically, for any \( x \in \mathbb{R}^{p_1 \times \cdots \times p_m} \), we have
\[
x = \sum_{k=1}^{r} \lambda \lo k (x)  U_k^{(1,0)}(x)  \otimes U_k^{(2,0)}(x)  \otimes \cdots \otimes U_k^{(m,0)}(x)   + \epsilon,
\]
where  \( U_k^{(\ell,0)}(x) \in \mathbb{R}^{p_\ell} \) are the mode-\( \ell \) component vectors, 
and \( \epsilon \in \mathbb{R}^{p_1 \times \cdots \times p_m} \) is a residual tensor. Similar to   the case of matrix-valued $x$, we assume that the first components of $U \lo k \hi {(1,0)}, \ldots, U \lo k \hi {(r-1,0)}$ are positive, and the first component of $U \lo k \hi {(r,0)}$ is such that $\lambda \lo k (x)$ is positive. Let $U \lo k \hi {(\ell)} = \lambda \lo k \hi {1/r} (x) U \lo k \hi {(\ell,0)}$. Then 
\[
x = \sum_{k=1}^{r}   U_k^{(1)}(x)  \otimes U_k^{(2)}(x)  \otimes \cdots \otimes U_k^{(m)}(x)   + \epsilon.
\]
Analogous to the development in Section~4, 
let \( \mathcal{H}_1, \cdots, \mathcal{H}_m \) denote RKHSs induced by the kernels \( \kappa_1, \cdots, \kappa_m \), respectively, 
and define the feature map 
\( F : \Omega_X \to \mathcal{H}_1 \otimes \cdots \otimes \mathcal{H}_m \):
\[
F (x)  = \sum_{k=1}^r \kappa_1(\cdot,U_k^{(1)}(x)) \cdots  \kappa_m(\cdot,U_k^{(m)}(x) ) \in  \mathcal{H}_1 \otimes \cdots \otimes \mathcal{H}_m.
\]
With the feature map thus defined, the covariance operators \(\Sigma \lo {FF}: \ca H \lo 1 \otimes \cdots \otimes \ca H \lo m  \to   \ca H \lo 1 \otimes \cdots \otimes \ca H \lo m\), \(\Sigma \lo {FY}: \ca H \lo Y \to  \ca H \lo 1 \otimes \cdots \otimes \ca H \lo m\) and \(R \lo {FY}: \ca H \lo Y \to  \ca H \lo 1 \otimes \cdots \otimes \ca H \lo m\) 
are defined in the same way as in the matrix $X$ case, by (\ref{equ:Covariance}) and Definition~\ref{definition:regression}. 
We define 
\(s(Y) = \Sigma_{FF}^{-1} \Sigma_{FY} (\kappa_Y(\cdot, Y) - E [\ka \lo Y (\cdot, Y)])\), which is a member of 
\(\mathcal{H}_1 \otimes \cdots \otimes \mathcal{H}_m\). To facilitate mode-wise iterative minimization, for $k = 1, \ldots, m$, we define the linear operators,
\(T_{s(Y)}^{(k)} : \ca H \lo {-k} \to  \mathcal{H}_k, \)
where 
$ 
\ca H \lo {-k}$ is the space $ \mathcal{H}_1 \otimes \cdots \otimes \mathcal{H}_{k-1} \otimes \mathcal{H}_{k+1} \otimes \cdots \otimes \mathcal{H}_m$ by the bounded bilinear forms
\begin{align*}
T \lo {s(Y)} \hii k: (  \ca H \lo k, \ca H \lo {-k} ) \to \real, \quad (f \hii k, f \hii {-k}) \mapsto \langle f \hii k f \hii {-k}, s(Y) \rangle \lo {\ca H \lo k \otimes \ca H \lo {-k}}, \quad  k = 1, \ldots, m,  
\end{align*}
where $f \hii k \in \ca H \lo k$ and $f \hii {-k} \in \ca H \lo {-k}$. We seek to minimize the least-squares criterion 
\begin{align*}
E    \left(\| s(Y) - u\lo 1 (Y) f \lo 1 \hii 1 \cdots f \lo 1 \hii m - \cdots - 
     u\lo d (Y) f \lo d \hii 1 \cdots f \lo d \hii m \| \lo {\ca H \lo 1 \otimes \cdots \otimes \ca H \lo m} \hi 2 \right) \end{align*}
subject to the constraints that \(\{f \lo 1 \hii 1, \ldots, f \lo d \hii 1 \}, \cdots,     \{f \lo 1 \hii m, \ldots, f \lo d \hii m \}\) are orthonormal sets in $\ca H \lo 1, \ldots, \ca H \lo m$, respectively. The optimization can be performed iteratively, 
with each step admitting an explicit solution. 
For instance, the optimal 
$f_1^{(1)*}, \ldots, f_1^{(m)*}$ 
can be obtained by the following corollary, 
which parallels Corollary~\ref{corollary:explicit sol for cp}. 
In what follows, let 
$f_1^{(-k)}$ denote the $(m-1)$-dimensional vector of functions 
\((f_1^{(1)}, \ldots, f_1^{(k-1)}, f_1^{(k+1)}, \ldots, f_1^{(m)})\). 
Define
\begin{align*}
    L_C(f^{(1)}, \ldots, f^{(m)}, u)
    = E\big(\| s(Y) - u(Y) f^{(1)} \cdots f^{(m)} \|_{\otimes}^2\big),
\end{align*}
where the subscript \(C\) in \(L_C\) indicates the “CP-approach.”

\begin{corollary}\label{corollary:explicit sol for cp high-order}  \
\begin{enumerate}
   \item For each $k \in \{1, \ldots, m \}$, each   $u \in L \lo 2 (P \lo Y)$,  and each  $f   \hii {-k} \in \ca H \lo {-k}$, the minimizer of $L \lo C(f \hii 1, \ldots, f \hii m,u)$ over $f \hii k \in \ca H \lo k $ is 
\begin{align*}
    f \hi {(k)*}  = E  [u (Y)  T \lo {s(Y)} \hi k  f \hii {-k}   ] /    E  [u \hi 2 (Y) \| f \hii {-k} \| \hi 2 \lo {\ca H \lo {-k}} ]. 
\end{align*}
\item For given $f \hii 1 \in \ca H \lo 1, \ldots, f \hii m \in \ca H \lo m $, the minimizer of $L \lo C (f \hii 1, \ldots, f \hii m ,u)$ over $u \in L \lo 2 (P \lo Y)$  is 
\begin{align*} 
   u \hi * (y)  = 
    \langle f \hii 1 \cdots f \hii m , s(y) \rangle \lo   {\ca H \lo 1 \otimes \cdots \otimes \ca H \lo m}/ \|f \hii 1 \cdots f \hii m  \| \lo {\ca H \lo 1 \otimes \cdots \otimes \ca H \lo m} \hi 2 . 
\end{align*}
\end{enumerate}
\end{corollary}
The proof of this corollary follows directly from the argument of 
Corollary~\ref{corollary:explicit sol for cp} and is therefore omitted. 
From this point onward, the development of the sample-level implementation—
including orthogonalization, GCV, 
and iterative minimization—is entirely analogous to the matrix-valued 
\( X \) case and is thus omitted for brevity.

\section{Simulation Studies} We consider three simulation settings for  Tucker-form NTSDR and one simulation setting for the CP-form NTSDR. 
For the first three Tucker-form NTSDR settings, our response variable $Y$ is set to be 
\begin{align*}
Y   = \log \left[1+\tsum_{k=1}^{s}\tsum_{l=1}^{t}\tsum_{i=1}^r f_k(U_i) g_l(V_i) \right]+\varepsilon, 
\end{align*}
where $\{(U \lo i, V \lo i ): i = 1, \ldots, r\}$ are obtained from the singular value decomposition of $X$, as described in Section 2. 
The three simulation settings differ in the choices of ${(f_k, g_\ell)}$ as follows:
\begin{enumerate}
    \item[\textbf{I.}] \( s = t = 1 \): 
    $f_1(u) = u_1^3$, and $g_1(v) = v_1^5$.
    \item[\textbf{II.}] \( s = t = 1 \):    
    $f_1(u) = u_1^3 \exp(1 + u_2)$, and $g_1(v) = v_1^5 \exp(1 + v_2)$.    
    \item[\textbf{III.}] \( s = t = 2 \):   $f_1(u) = u_1^3$,$ f_2(u) = \exp(1 + u_2)$,  $g_1(v) = v_1^5$, and $g_2(v) = \exp(1 + v_2)$.  
\end{enumerate}
For the fourth simulation setting for the CP-form NTSDR, we set the response to be 
\begin{align*}
 Y = \log \left[ 1 + \tsum_{k=1}^{d} \tsum_{i=1}^r f_k(U_i) g_k(V_i) \right] 
   \end{align*}
   where 
\begin{enumerate}
    \item[\textbf{IV.}] $d = 4$, 
    $f_1(u) = u_1^3$, 
    $f_2(u) = \exp(1 + u_2)$, 
    $f_3(u) = \log(1 + u_1)$, 
    $f_4(u) = u_1 + u_2$, 
    $g_1(v) = v_1^5$, 
    $g_2(v) = \exp(1 + v_2)$, 
    $g_3(v) = \log(1 + v_1)$, 
    and $g_4(v) = v_1 + v_2$.
\end{enumerate}

In all settings, the noise term $\varepsilon$ is Gaussian with mean 0 and variance determined by the signal-to-noise ratio  $\var (Y) / \var (\varepsilon) = 0.25$. In practice, we first compute the empirical variance of the noiseless signal term and then scale $\sigma^2$ so that the resulting ratio equals $0.25$. 
This calibration ensures that the signal variance is one quarter of the noise variance, providing a moderately noisy yet informative setting for evaluating the proposed methods. The predictor matrix \(X\) is generated under the following three scenarios: 
\begin{enumerate}
    \item[(A)] independent Gaussian predictors, $N(\Theta, I_{p \times q})$,
    \item[(B)] independent non-Gaussian predictors, 
    \(\tfrac{1}{2} N(\Theta - \mathbf{1}_{p \times q}, 2 I_{p \times q}) + \tfrac{1}{2} N(\Theta + \mathbf{1}_{p \times q}, 2 I_{p \times q}),\)
    \item[(C)] correlated Gaussian predictors, $N(\Theta, 0.8 I_{p \times q} + 0.2 \mathbf{1}_{p \times q} \mathbf{1}\trans_{p \times q})$. 
\end{enumerate}
We take \((p, q) = (5,5), (10,10)\),   \(n = 100, 200\). Also, \(\Theta\) is defined as follows: in \(\mathbf{I}\), \(\Theta_{11} = 5\) and all other entries are zero; in \(\mathbf{II}\), \(\mathbf{III}\), and \(\mathbf{IV}\),  \(\Theta_{11} = \Theta_{22} = 5\), and the rest are zero.

To evaluate prediction accuracy, we compute the fitted response 
$\hat{Y}$ by plugging the estimated nonlinear functions 
$\{\hat f_k, \hat g_\ell\}$ into the population model:
\[
\hat{Y}
=
\log\!\left[
1+
\sum_{i=1}^{r}
\Big(
\sum_{k=1}^{s}\sum_{\ell=1}^{t}\hat f_k(U_i)\hat g_\ell(V_i)
\;\;\text{or}\;\;
\sum_{k=1}^{d}\hat f_k(U_i)\hat g_k(V_i)
\Big)
\right],
\]
for the Tucker- and CP-form NTSDR models, respectively.
The estimated functions $\hat f_k$ and $\hat g_\ell$ are obtained from the corresponding sample-level optimization procedures. To assess the closeness between the estimated \(\hat{Y}\) and the true \(Y\) , we use distance correlation \citep{szekely2007measuring}, which quantifies the strength of statistical dependence between vectors or matrices. We compare our methods (Tucker- and CP-forms of NTSDR) with GSIR \citep{Lee2013}  the Kernel SIR (KSIR) \citep{yeh2008}. In Table 1 we report the means and standard errors of the distance correlation between $\hat Y$ and $Y$ for the four methods evaluated from test data samples of size $n \lo t = 100$, and based on 200 replications. As shown in Table~\ref{tab:tab1}, our methods consistently outperform GSIR and KSIR, with comparable performance between the two versions of the NTSDR.

\begin{table}[htbp]
\centering
\caption{Distance correlation between true and estimated responses across four settings and three predictor designs.}
\label{tab:tab1}
\begin{adjustbox}{max width=\textwidth}
\begin{tabular}{c|c|cccc|cccc}
\toprule
& & \multicolumn{8}{c}{Different samples and dimensions $(n,p,q)$}\\
\midrule
& & \multicolumn{4}{c|}{$(100,5,5)$} & \multicolumn{4}{c}{$(200,5,5)$}\\
\midrule
 & X & NTSDR-Tu & NTSDR-CP & GSIR & KSIR & NTSDR-Tu & NTSDR-CP & GSIR & KSIR \\
\midrule
I & A & \multicolumn{2}{c}{0.771(0.076)} & 0.609(0.109) & 0.649(0.095) & \multicolumn{2}{c}{0.742(0.183)} & 0.743(0.177) & 0.663(0.102)\\
  & B & \multicolumn{2}{c}{0.716(0.088)} & 0.704(0.097) & 0.682(0.094) & \multicolumn{2}{c}{0.797(0.060)} & 0.744(0.093) & 0.652(0.065) \\
  & C & \multicolumn{2}{c}{0.693(0.064)} & 0.661(0.099) & 0.549(0.085) & \multicolumn{2}{c}{0.740(0.086)} & 0.703(0.094) & 0.627(0.071)\\
II & A & \multicolumn{2}{c}{0.637(0.150)} & 0.359(0.141) & 0.313(0.132) & \multicolumn{2}{c}{0.614(0.137)} & 0.386(0.134) & 0.407(0.129)\\
& B & \multicolumn{2}{c}{0.642(0.168)} & 0.325(0.118) & 0.340(0.127) & \multicolumn{2}{c}{0.735(0.141)} & 0.369(0.139) & 0.400(0.153)\\
& C & \multicolumn{2}{c}{0.617(0.139)} & 0.327(0.111) & 0.289(0.127) & \multicolumn{2}{c}{0.673(0.158)} & 0.349(0.124) & 0.321(0.143)\\
III & A & 0.851(0.069) & 0.853(0.072) & 0.398(0.137) & 0.433(0.127) & 0.870(0.074) & 0.873(0.057) & 0.430(0.163) & 0.441(0.165)\\
& B & 0.876(0.080) & 0.822(0.090) & 0.404(0.158) & 0.381(0.142) & 0.880(0.058) & 0.866(0.060) & 0.462(0.143) & 0.459(0.160)\\
& C & 0.835(0.075) & 0.862(0.066) & 0.419(0.158) & 0.394(0.135) & 0.856(0.052) & 0.858(0.067) & 0.454(0.160) & 0.422(0.179)\\
IV & A & 0.844(0.082) & 0.796(0.092) & 0.407(0.157) & 0.396(0.138) & 0.836(0.091) & 0.803(0.144) & 0.437(0.166) & 0.423(0.154)\\
& B & 0.868(0.064) & 0.789(0.092) & 0.213(0.100) & 0.188(0.076) & 0.858(0.071) & 0.809(0.058) & 0.448(0.171) & 0.415(0.148)\\
& C & 0.808(0.086) & 0.791(0.135) & 0.448(0.157) & 0.400(0.137) & 0.829(0.118) & 0.773(0.172) & 0.471(0.163) & 0.446(0.166)\\
\midrule
& & \multicolumn{4}{c|}{$(100,10,10)$} & \multicolumn{4}{c}{$(200,10,10)$}\\
\midrule
 & X & NTSDR-Tu & NTSDR-CP & GSIR & KSIR & NTSDR-Tu & NTSDR-CP & GSIR & KSIR \\
\midrule
I & A & \multicolumn{2}{c}{0.766(0.082)} & 0.462(0.104) & 0.195(0.058) & \multicolumn{2}{c}{0.764(0.096)} & 0.567(0.094) & 0.502(0.098)\\
& B &  \multicolumn{2}{c}{0.762(0.083)} & 0.487(0.085) & 0.204(0.059) &  \multicolumn{2}{c}{0.782(0.079)} & 0.520(0.100) & 0.494(0.093) \\
& C &  \multicolumn{2}{c}{0.609(0.201)} & 0.398(0.095) & 0.194(0.052) &  \multicolumn{2}{c}{0.705(0.178)} & 0.455(0.114) & 0.371(0.089) \\
II & A &  \multicolumn{2}{c}{0.595(0.172)} & 0.304(0.094) & 0.198(0.060) &  \multicolumn{2}{c}{0.645(0.166)} & 0.382(0.099) & 0.303(0.095)\\
& B &  \multicolumn{2}{c}{0.576(0.189)} & 0.285(0.091) & 0.178(0.057) &  \multicolumn{2}{c}{0.632(0.177)} & 0.381(0.094) & 0.322(0.106)\\
& C &  \multicolumn{2}{c}{0.556(0.179)} & 0.265(0.085) & 0.195(0.051) &  \multicolumn{2}{c}{0.602(0.171)} & 0.304(0.108) & 0.258(0.097)   \\
III & A & 0.830(0.061) & 0.846(0.076) & 0.226(0.072) & 0.190(0.085) & 0.838(0.058) & 0.816(0.066) & 0.252(0.097) & 0.225(0.102)\\
& B & 0.835(0.060) & 0.839(0.066) & 0.229(0.080) & 0.194(0.080) & 0.839(0.060) & 0.829(0.060) & 0.268(0.100) & 0.221(0.101)  \\
& C & 0.832(0.058) & 0.830(0.064) & 0.259(0.097) & 0.180(0.090) & 0.843(0.061) & 0.819(0.061) & 0.250(0.101) & 0.243(0.103)\\
IV & A &0.828(0.055) & 0.783(0.103) & 0.242(0.088) & 0.190(0.075) & 0.853(0.059) & 0.797(0.057) & 0.264(0.100) & 0.230(0.097)\\
& B & 0.822(0.061) & 0.799(0.087) & 0.250(0.089) & 0.177(0.077) & 0.846(0.062) & 0.794(0.060) & 0.271(0.100) & 0.238(0.102)  \\
& C & 0.869(0.058) & 0.820(0.078) & 0.271(0.092) & 0.179(0.068) & 0.872(0.062) & 0.822(0.057) & 0.278(0.106) & 0.234(0.067)\\
\bottomrule
\end{tabular}
\end{adjustbox}
\end{table}

To further examine the performance of our methods in recovering the true sufficient predictor, we also computed the distance correlations between the true and estimated sufficient predictors. 
Specifically, for each observation, we compute the true interaction features 
\[
S_{\mathrm{true}} = \bigl[f_{k}(u_{ij})\, g_{\ell}(v_{ij})\bigr]_{i,j,k,\ell}
\quad\text{and}\quad
S_{\mathrm{est}} = \bigl[\hat f_{k}(u_{ij})\, \hat g_{\ell}(v_{ij})\bigr]_{i,j,k,\ell},
\]
each reshaped into a matrix of size $nr \times (st)$ for the Tucker-form model (or $nr \times d$ for the CP-form model). 
We then compute the distance correlation between the columns of $S_{\mathrm{true}}$ and $S_{\mathrm{est}}$ as a global measure of the accuracy of sufficient predictor recovery: 
\[
\mathrm{dCor}(S_{\mathrm{true}}, S_{\mathrm{est}}) 
= \frac{\mathrm{dCov}(S_{\mathrm{true}}, S_{\mathrm{est}})}
{\sqrt{\mathrm{dVar}(S_{\mathrm{true}})\,\mathrm{dVar}(S_{\mathrm{est}})}},
\]
where $\mathrm{dCov}$ and $\mathrm{dVar}$ denote the distance covariance and variance as defined in \cite{szekely2007measuring}. 
This criterion captures both linear and nonlinear associations between the recovered and true {sufficient predictors}. 
Table~\ref{tab:tab2} reports the mean and standard errors of the resulting distance correlations evaluated from test data samples of size $n \lo t = 100$, and based on 200 replications.
The results demonstrate that both Tucker- and CP-form NTSDR methods successfully recover the intrinsic functional structure across all simulation scenarios, thus confirming the fidelity of the recovered nonlinear components. These findings collectively indicate that NTSDR effectively captures both global response variation and local structural dependencies in the underlying tensor predictors. In addition to the regression settings discussed above, we also conduct a simulation study for a classification setting, which is presented in the supplement.

\begin{table}[htbp]
    \centering
    \caption{Distance correlation between true and estimated structures across four settings and three predictor designs.}
    \label{tab:tab2}
    \begin{adjustbox}{max width=\textwidth}
    \begin{tabular}{l l c c c c c c c c}
        \toprule
        & & \multicolumn{8}{c}{Different samples and dimensions $(n,p,q)$}\\
        \cmidrule(lr){3-10}
         & X & \multicolumn{2}{c}{$(100,5,5)$} & \multicolumn{2}{c}{$(200,5,5)$} & \multicolumn{2}{c}{$(100,10,10)$} & \multicolumn{2}{c}{$(200,10,10)$}\\
        \midrule
        & & \multicolumn{2}{c}{NTSDR-Tu \& CP} & \multicolumn{2}{c}{NTSDR-Tu \& CP} & \multicolumn{2}{c}{NTSDR-Tu \& CP} & \multicolumn{2}{c}{NTSDR-Tu \& CP}\\
        \midrule
        I & A & \multicolumn{2}{c}{0.981 (0.047)} & \multicolumn{2}{c}{0.942 (0.049)} & \multicolumn{2}{c}{0.906 (0.066)} & \multicolumn{2}{c}{0.932 (0.068)}\\
          & B & \multicolumn{2}{c}{0.936 (0.060)} & \multicolumn{2}{c}{0.983 (0.055)} & \multicolumn{2}{c}{0.894 (0.088)} & \multicolumn{2}{c}{0.939 (0.082)}\\
          & C & \multicolumn{2}{c}{0.981 (0.061)} & \multicolumn{2}{c}{0.974 (0.061)} & \multicolumn{2}{c}{0.735 (0.229)} & \multicolumn{2}{c}{0.744 (0.205)}\\
        II & A & \multicolumn{2}{c}{0.644 (0.190)} & \multicolumn{2}{c}{0.706 (0.173)} & \multicolumn{2}{c}{0.633 (0.186)} & \multicolumn{2}{c}{0.682 (0.181)}\\
          & B & \multicolumn{2}{c}{0.693 (0.215)} & \multicolumn{2}{c}{0.802 (0.182)} & \multicolumn{2}{c}{0.616 (0.208)} & \multicolumn{2}{c}{0.687 (0.207)}\\
          & C & \multicolumn{2}{c}{0.688 (0.187)} & \multicolumn{2}{c}{0.662 (0.224)} & \multicolumn{2}{c}{0.553 (0.226)} & \multicolumn{2}{c}{0.594 (0.197)}\\
        \midrule
         & & NTSDR-Tu & NTSDR-CP & NTSDR-Tu & NTSDR-CP & NTSDR-Tu & NTSDR-CP & NTSDR-Tu & NTSDR-CP\\
        \midrule
        III & A & 0.855(0.122) & 0.858(0.116) & 0.850(0.119) & 0.844(0.118) & 0.842(0.095) & 0.836(0.091) & 0.817(0.088) & 0.844(0.089)\\
          & B & 0.827(0.136) & 0.838(0.119) & 0.875(0.135) & 0.828(0.119) & 0.843(0.104) & 0.827(0.107) & 0.831(0.098) & 0.846(0.101)\\
          & C & 0.818(0.147) & 0.791(0.112) & 0.821(0.135) & 0.797(0.104) & 0.792(0.112) & 0.835(0.094) & 0.832(0.096) & 0.836(0.086)\\
        IV & A & 0.821(0.133) & 0.840(0.116) & 0.815(0.121) & 0.838(0.123) & 0.834(0.109) & 0.794(0.117) & 0.815(0.107) & 0.809(0.117)\\
          & B & 0.812(0.135) & 0.815(0.110) & 0.806(0.140) & 0.862(0.120) & 0.802(0.113) & 0.824(0.119) & 0.829(0.107) & 0.799(0.125)\\
          & C & 0.772(0.149) & 0.825(0.112) & 0.817(0.138) & 0.787(0.127) & 0.848(0.104) & 0.739(0.134) & 0.808(0.092) & 0.786(0.111)\\
        \bottomrule
    \end{tabular}
    \end{adjustbox}
\end{table}

\section{Real Data Applications}
We have applied the proposed NTSDR methods to two datasets, one for the regression setting and another for the classification setting, and both with matrix-valued predictors. The first application, presented in this section, is a Categorical Subjective Image Quality (CSIQ) dataset; the second, presented in the Supplementary Material, is an electroencephalography (EEG) dataset. 
We also compared our proposed methods with several existing dimension reduction techniques through these applications. 

The CSIQ  dataset, available at \url{https://s2.smu.edu/~eclarson/csiq.html}, is a widely used benchmark for evaluating image quality models.  It includes 30 reference images, each degraded by six types of distortions (e.g., Gaussian noise, JPEG compression, contrast change) at varying severity levels, yielding approximately 28–29 distorted versions per reference. Subjective quality ratings were collected under controlled viewing conditions and quantified using Differential Mean Opinion Scores (DMOS). In this study, we focus on grayscale images to isolate the impact of distortions from color effects. Since each image is represented as a  \( 512 \times 512 \) matrix, direct analysis is computationally intensive. To mitigate this, we apply an initial dimensionality reduction step by projecting each image into a lower-dimensional space \(\mathbb{R}^{d \times d}\) , where  \( d \)  is chosen to balance computational efficiency with information preservation.

To evaluate generalization, we partition the dataset into training and testing sets. For each reference image, 22 of the 29 distortion-severity combinations are used for training, and the others for testing. We compute the Pearson correlation between the estimated quality scores  \(\hat{Y} \)  and the ground-truth DMOS values  \(Y \), repeating this process over 10 trials. Table~\ref{table:CSIQ4} compares one of the proposed method CP-NTSDR with GSIR. CP-NTSDR outperforms GSIR on 24 of the 30 reference images and also shows comparable results in three of the remaining cases. The consistently high correlation values indicate strong generalization to unseen distortions. The superior performance of NTSDR-CP can be attributed to its ability to preserve intrinsic data structures while substantially reducing the number of model parameters. Unlike GSIR, which applies a global transformation to the entire predictor space, our method identifies a compact and structured subspace that retains the dominant variations in perceptual quality.
\begin{table}[ht!]
\centering
\scriptsize
\caption{Mean Pearson correlation between estimated image quality scores and ground-truth DMOS values.}
\label{table:CSIQ4}
\begin{tabular}{lcc|lcc}
\toprule
Image & NTSDR-CP & GSIR & Image & NTSDR-CP & GSIR \\
\midrule
1600 & \textbf{0.875} & 0.824 & Aerial\_city & \textbf{0.765} & 0.670 \\
Boston & 0.726 & \textbf{0.721} & Bridge & \textbf{0.707} & 0.673 \\
Butter\_flower & \textbf{0.781} & 0.707 & Cactus & \textbf{0.733} & 0.692 \\
Child\_swimming & \textbf{0.827} & 0.722 & Couple & 0.495 & \textbf{0.675} \\
Elk & \textbf{0.844} & 0.783 & Family & \textbf{0.818} & 0.762 \\
Fisher & \textbf{0.901} & 0.852 & Foxy & \textbf{0.837} & 0.755 \\
Geckos & \textbf{0.745} & 0.665 & Lady\_liberty & 0.573 & \textbf{0.663} \\
Lake & \textbf{0.613} & 0.470 & Log\_seaside & \textbf{0.796} & 0.729 \\
Monument & \textbf{0.792} & 0.698 & Native\_american & \textbf{0.774} & 0.775 \\
Redwood & \textbf{0.744} & 0.679 & Roping & \textbf{0.819} & 0.792 \\
Rushmore & \textbf{0.557} & 0.516 & Shroom & 0.847 & \textbf{0.853} \\
Snow\_leaves & \textbf{0.804} & 0.548 & Sunset\_sparrow & \textbf{0.935} & 0.786 \\
Sunsetcolor & \textbf{0.932} & 0.788 & Swarm & \textbf{0.852} & 0.794 \\
Trolley & \textbf{0.809} & 0.780 & Turtle & 0.842 & \textbf{0.848} \\
Veggies & \textbf{0.739} & 0.703 & Woman & \textbf{0.809} & 0.779 \\
\bottomrule

\end{tabular}
\end{table}
These findings underscore the efficacy of the NTSDR methods in analyzing tensor-structured data. By preserving matrix or tensor structures throughout the modeling process, the NTSDR yields an interpretable and computationally efficient framework for uncovering meaningful low-dimensional representations.
 
\section{Conclusion}
We have developed a unified framework for nonlinear SDR of matrix- and tensor-valued predictors that preserves the tensor structure of the data. This is achieved by constructing a tensor regression operator that maps from the response RKHS to the tensor feature space, whose range space is guaranteed to generate the same $\sigma$-field as that generated by the tensor-preserving sufficient predictors.  Our results extend the existing tensor-preserving linear  SDR methods by following two well-known approaches to tensor decomposition: the Tucker decomposition and the CP decomposition.  We developed the population-level Fisher consistency, convergence rates, and asymptotic consistency, as well as implementation algorithms for estimation and tuning.  Besides preserving the tensor interpretation, our methods offer significant improvement in estimation accuracy, as they substantially reduce the number of parameters involved, which is shown in both simulation studies and real applications.

\bibliographystyle{agsm}
\bibliography{ref}

\end{document}